\documentclass{amsart}
\usepackage{amssymb,latexsym,graphicx,psfrag,epstopdf,hyperref}

\usepackage[all,arc,dvips]{xy}

\newtheorem{proposition}{Proposition}[section]
\newtheorem{theorem}[proposition]{Theorem}

\newtheorem{question}[proposition]{Question}

\newtheorem{corollary}[proposition]{Corollary}
\newtheorem{conjecture}[proposition]{Conjecture}
\newtheorem{lemma}[proposition]{Lemma}
\newtheorem*{acknowledgments}{Acknowledgments}
\def\frac#1#2{{#1\over#2}}

\numberwithin{equation}{section}

\def\Area{\mbox{\rm{Area}}}
\def\Vol{\mbox{\rm{Vol}}}
\def\PSL{\mbox{\rm{PSL}}}
\def\Isom{\mbox{\rm{Isom}}}
\def\Stab{\mbox{\rm{Stab}}}

\def\co{\colon\thinspace}

\def\cal#1{\mathcal{#1}}
\def\frak#1{\mathfrak{#1}}
\def\bb#1{\mathbb{#1}}
\def\th{\theta}

\def\p{\pi}

\def\C{\Gamma}
\def\D{\Delta}

\begin{document}
\Large
\title{Immersed Turnovers In Hyperbolic $3$--Orbifolds}
\author{Shawn~Rafalski}
\address{Department of Mathematics and Computer Science, Fairfield University, Fairfield, CT 06824}
\email{srafalski@fairfield.edu}
\date{December 2008}

\begin{abstract}
We show that any immersion, which is not a covering of an embedded $2$--orbifold, of a totally geodesic hyperbolic turnover in a complete orientable hyperbolic $3$--orbifold  is contained in a hyperbolic $3$--suborbifold with totally geodesic boundary, called the ``turnover core,'' whose volume is bounded from above by a function depending only on the area of the given turnover.  Furthermore, we show that, for a given type of turnover, there are only finitely many possibilities for the turnover core. As a corollary, if the volume of a complete orientable hyperbolic $3$--orbifold is at least $2 \pi$ and if the fundamental group of the orbifold contains the fundamental group of a hyperbolic turnover (i.e., a triangle group), then the orbifold contains an embedded hyperbolic turnover.
\end{abstract}

\maketitle

\section{Introduction}\label{S:Intro}
It is well known that the thrice-punctured sphere is the only orientable hyperbolic surface which is rigid, in the sense that it admits a unique complete hyperbolic structure. The analogue of this surface in the orbifold setting is the hyperbolic turnover. A \emph{hyperbolic turnover} is the double, along the boundary, of a hyperbolic triangle whose interior angles are integer submultiples of $\pi$. As an orbifold, a hyperbolic turnover is topologically a $2$--sphere with three cone points whose orders correspond to the submultiples of $\pi$ in the associated hyperbolic triangle.  Like the thrice-punctured sphere, hyperbolic turnovers admit a unique complete hyperbolic structure. The goal of this paper is to prove the following theorem, which can be viewed either as a finiteness result or as the turnover analogue of several well-known theorems from classical 3--manifold topology (see Corollaries \ref{C:TorusThmEsque} and \ref{C:JSJDecompEsque} below).

\begin{theorem}\label{T:ImmersedTurnover}
	Let $f\co \cal{T} \to Q$ be a totally geodesic (equivalently, $\pi_{1}$--injective) immersion of a compact, hyperbolic turnover $\cal{T} = \cal{T}(p,q,r)$ in an orientable hyperbolic $3$--orbifold $Q$. Assume that $f(\cal{T})$  does not cover an embedded turnover or an embedded triangle with mirrored sides. Then $Q$ contains a finite (possibly empty) collection $\{ \cal{T}_{i} \}$ of embedded, pairwise disjoint, totally geodesic hyperbolic turnovers (and totally geodesic hyperbolic triangles with mirrored sides)  satisfying the following: \smallskip
	\begin{enumerate}
		\item $f(\cal{T}) \cap \cal{T}_{i} = \emptyset$ for each $i$
		\item The number of turnovers (and triangles with mirrored sides) in the collection 
			$\{ \cal{T}_{i} \}$ is bounded above by a function of $p,q$ and $r$
		\item If $n$ is the order of a cone point of any turnover in the collection $\{ \cal{T}_{i} \}$, 
			then $$n \in \left\{ 2,3,...,9,p,q,r,2p,2q,2r \right\}.$$
		\item If $\p / n$ is an angle of a triangle with mirrored sides in the 
			collection $\{ \cal{T}_{i} \}$, then $$n \in \left\{ 2,3,...,9,p,q,r,2p,2q,2r \right\}.$$
		\item  If $Q'$ is the component of $Q -  \cup_{i} \cal{T}_{i}$ which contains $f(\cal{T})$, then the 
			metric closure of $Q'$ is a small hyperbolic $3$--orbifold with (possibly empty) totally 
			geodesic boundary. If $\partial Q'$ is not empty, then $\Vol(Q') < H \cdot  \Area(\cal{T}),$
			where $H=1.199678...$ is the positive solution of the equation $x = \coth x.$ 
			If $\partial Q'$ is empty, then $\Vol(Q') < \Area(\cal{T})$.
	\end{enumerate}
	Furthermore, for a given $(p,q,r)$--turnover, there are only finitely many possibilities for the orbifold $Q'$ described above.
\end{theorem}  

Because the area of any hyperbolic turnover is bounded above by $2 \p$, we have the following

\begin{corollary}\label{C:TorusThmEsque}
	(The Turnover Theorem) Let $Q$ be an orientable hyperbolic $3$--orbifold with $\Vol(Q) \geq 2 \p$. Then $Q$ contains an embedded turnover (or an embedded triangle with mirrored sides), if it contains an immersed turnover.
\end{corollary}

Dunbar \cite{Dunbar88-1} showed that every compact, irreducible, atoroidal 3--orbifold can be split (uniquely, up to isotopy) along a system of essential, pairwise non-parallel hyperbolic turnovers into pieces which contain no essential (embedded) turnovers. The next corollary, which follows from the arguments of Section \ref{S:NbhdHyp}, says that the Dunbar decomposition is the turnover analogue of the JSJ-decomposition of a 3--manifold.

\begin{corollary}\label{C:JSJDecompEsque}
	Let $Q$ be a compact, irreducible, orientable, atoroidal $3$--orbifold. Then any immersion $f\co \cal{T} \to Q$ of a hyperbolic turnover into $Q$ is homotopic into a unique component of the Dunbar decomposition, up to parallel boundary components of the decomposition.
\end{corollary}

The motivation for Theorem \ref{T:ImmersedTurnover} begins with a question which Kirby attributes to Martin \cite[Problem $3.70$]{Kirby96Problems}:

\begin{question}\label{Q:MartinConj}
Given a Kleinian group $\C$ and a turnover subgroup $T$ with invariant plane $\Pi$ in $\bb{H}^{3}$, is it true that for all $\gamma \in \C$ either $\gamma \Pi = \Pi$ or $\gamma \overline{\Pi} \cap \overline{\Pi} = \emptyset$? Equivalently, is it true that the turnover $\cal{T} = \Pi / T$  covers an embedded turnover in the orbifold $\bb{H}^{3} / \C$?
\end{question}

\begin{figure}[htbp]
		\begin{center}
		\psfrag{A}{{\small $A$}}
		\psfrag{B}{{\small $B$}}
		\psfrag{C}{{\small $C$}}
		\psfrag{D}{{\small $D$}}
		\psfrag{2}{{\small $2$}}
		\psfrag{3}{{\small $3$}}
		\psfrag{4}{{\small $4$}}
		\psfrag{5}{{\small $5$}}
		\psfrag{pi/4}{{\small ${\displaystyle \p \over \displaystyle 4}$}}
		\psfrag{pi/5}{{\small ${\displaystyle \p \over \displaystyle 5}$}}
		\scalebox{1}{\includegraphics{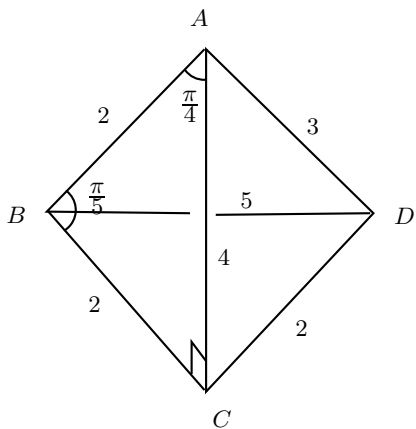}}
		\caption{A compact hyperbolic tetrahedron which yields a $3$--orbifold with an immersed 
			turnover}
		\label{F:BaskMacbIntro}
		\end{center}
	\end{figure}

As it turns out, the answer to this question is no, although the author believes that this must have been known to Martin at the time that Kirby added this question to his list of problems in low dimensional topology. Consider the tetrahedron in Figure \ref{F:BaskMacbIntro}.
The integers at the edges indicate the dihedral angles as submultiples of $\pi$, so, for instance, an edge labeled $2$ has a dihedral angle of $\pi / 2$.  It is known that this tetrahedron can be realized in hyperbolic $3$--space, and that the group generated by reflections in its faces is a discrete group of isometries.  Let $\C_{3}$ be the orientation-preserving subgroup of index two inside this reflection group.  Consider the face $ABC$. It is not difficult to see that this face is a hyperbolic triangle with the angles $\pi /2, \pi/4,$ and $\pi /5$, as indicated in the figure.  Consider the centralizer, in the reflection group, of reflection in this face.  Baskan and Macbeath \cite{BaskanMacbeath} proved that the orientation-preserving subgroup of this centralizer is a turnover subgroup of $\C_{3}$. (A \emph{turnover group} is the orbifold fundamental group of a turnover.) The corresponding geometric turnover is immersed in the quotient $3$--orbifold, because the order $3$ edge $AD$ meets the plane containing the face $ABC$ at an oblique angle.  

There are nine compact hyperbolic tetrahedra for which the group of isometries generated by edge rotations is discrete. In all nine cases, the corresponding hyperbolic $3$--orbifolds contain immersed turnovers.  The other examples are not as easily seen as in the Baskan and Macbeath example above, however, for this tetrahedron is the only one with a triangular face whose interior angles are integer submultiples of $\pi$.  Maclachlan has classified the immersed turnovers in these ``tetrahedral'' groups  for eight of the nine cases \cite{Maclachlan96-1}, although he makes an error in one of the cases which we will rectify (Proposition \ref{P:TOsInT9}). In any event, all of these counterexamples to Question \ref{Q:MartinConj} are ``small'' in several respects.  We will consider them in some detail in Section \ref{SS:Maclachlan}.  

On the other hand, Martin has proved the affirmative answer to his question for at least one class of hyperbolic turnover \cite{Martin96-1}. Specifically, he showed that $(2,3,p)$ turnovers, where $p \geq 7$, are always embedded in a hyperbolic $3$--orbifold.  In joint work with Gehring \cite{GehringMartin94-2}, this result was utilized as part of a general program for locating the minimal co-volume Kleinian group, which is conjectured to be the orientation-preserving index two subgroup of a particular Coxeter reflection group (in fact, the reflection group corresponds to one of the nine tetrahedral examples mentioned above).  

Knowing that a turnover is embedded in a hyperbolic $3$--orbifold can be used to find lower bounds for the volume of the $3$--orbifold, by considering embedded $\delta$--neighborhoods of the turnover (e.g. \cite[Theorem 1.11]{Martin96-1}). Furthermore, Kleinian groups which contain a turnover subgroup (or equivalently, hyperbolic $3$--orbifolds which contain hyperbolic turnovers) frequently turn up as candidates for the solutions to several types of extremal problems in hyperbolic geometry. For instance, recent work of Gehring and Martin \cite{GehringMartin05} attempts to determine the Margulis constant for $\bb{H}^{3}$, and the elements of one class of groups they consider (so called $(p,q,r)$--Kleinian groups) frequently contain turnover subgroups.  

Finally, because turnovers have the unique property of rigidity among orientable $2$--orbifolds, they act as shields to the effects of any deformation of the hyperbolic structure of the ambient hyperbolic $3$--orbifold.  Specifically, if a turnover $\cal{T}$ separates off a cusp $C$ of the $3$--orbifold $Q$ from the rest of $Q$, then any hyperbolic Dehn surgery performed on $C$ affects only the geometry of the component of $Q - \cal{T}$ which contains $C$. This is a specific instance of a phenomenon known as \emph{geometric isolation} \cite{NeumannReid93-1}, \cite{Calegari01-1}.  It follows from Corollary \ref{C:TorusThmEsque} that once the volume of a cusped hyperbolic $3$--orbifold $Q$ is large enough, then knowing its fundamental group contains a turnover subgroup is equivalent to knowing that there are pieces of $Q$ which are left alone by many deformations.  In light of all of these considerations, a complete classification of immersed turnovers in hyperbolic $3$--orbifolds would be quite useful.  It is the hope of the author that Theorem \ref{T:ImmersedTurnover} is a first step toward that classification. 

\begin{acknowledgments}
The author is grateful to his advisor Ian Agol for giving so freely of his time to discuss mathematics in all its forms. Special thanks also go to Marc Culler and Peter Shalen, who generously provided their time and suggestions for the work contained herein. Finally, the author thanks the referee for invaluable feedback.  
\end{acknowledgments}

\section{Definitions and Notation}\label{S:Defs}
An \emph{$n$--orbifold} $O$ is a metrizable topological space in which every point has a neighborhood which is diffeomorphic either to the quotient of $\bb{R}^{n}$ by a finite group action or to the quotient of $\bb{R}^{n-1} \times [0,\infty)$ by a finite group action. Points with neighborhoods modeled on the latter type of quotients make up the \emph{boundary} $\partial O$ of the $n$--orbifold, which is itself an $(n-1)$--orbifold. An orbifold is called \emph{geometric} if its interior is diffeomorphic to the quotient of a model geometric space by a discrete group of isometries. We will not define the term \emph{model geometric space}, but will rather point out that the $n$--dimensional sphere, Euclidean space, and hyperbolic space (denoted, respectively, by $\bb{S}^{n}$, $\bb{E}^{n}$, and $\bb{H}^{n}$) are the examples of model geometric spaces in which we will be interested.  There are several extensive references for and introductions to the definitions in this section \cite{BMP03-1}, \cite{CoopHodgKer00}.      

In all of what follows, $Q$  is a complete orientable hyperbolic $3$--orbifold. That is, $Q$ is the quotient $\bb{H}^{3}/\C$ of hyperbolic $3$--dimensional space by a Kleinian group $\C$. A \emph{Kleinian group} is  a discrete and non-elelmentary subgroup of $\PSL_{2}(\bb{C})$, where we identify $\PSL_{2}(\bb{C})$ with the group of orientation-preserving isometries of $\bb{H}^{3}$ via the Poincar\'e extension. Alternatively, $\C$ may be thought of as a subgroup of $\PSL_{2}(\bb{C})$ whose action on $\bb{H}^{3}$ is properly discontinuous (although not necessarily free). The action of a group $G$ on a space $X$ is \emph{properly discontinuous} if any compact subset $K \subset X$ is taken completely off of itself by all but only finitely many elements of $G$.  We will refer to $\C$ as the \emph{(orbifold) fundamental group of} $Q$, and sometimes denote it by $\p_{1}(Q)$.  Denote the covering projection by $\Phi\co \bb{H}^{3} \to Q$.

Let $2 \leq p \leq q \leq r$ be positive integers satisfying $\frac{1}{p} + \frac{1}{q} + \frac{1}{r} < 1$ (in generalizations of the terminology presented here it is allowed for $p,q,$ and $r$ to take on the value $\infty$). Then there is a hyperbolic triangle $\D = \D(p,q,r)$, unique up to isometry, with interior angles $\frac{\p}{p}$, $\frac{\p}{q}$ and $\frac{\p}{r}$. The group generated by reflections in the sides of $\D$ is a discrete subgroup of isometries of $2$--dimensional hyperbolic space $\bb{H}^{2}$. Let $T(p,q,r)$ be the unique normal subgroup of index two inside this reflection group which acts on $\bb{H}^{2}$ by orientation-preserving isometries. Then $T(p,q,r)$ is generated by rotations by $\frac{2\p}{p}$, $\frac{2\p}{q}$ and $\frac{2\p}{r}$ around the corresponding vertices of $\D$ (in fact, any pair of these rotations generates this group). Call this a \emph{turnover group} of Isom($\bb{H}^{2}$). (A remark: These are commonly referred to as \emph{triangle} groups, and while the author has a great deal of respect for tradition, he feels somewhat compelled to refer to these groups by the likeness of their associated $2$--dimensional geometric objects.)  A \emph{turnover subgroup}  of a Kleinian group $\C$ is a subgroup which is isomorphic to some $T(p,q,r)$. We will denote it by $T = T(p,q,r)$.

A turnover subgroup $T(p,q,r)$ of a Kleinian group $\C$ is generated by three elliptic elements $\gamma_{p}$, $\gamma_{q}$ and $\gamma_{r}$ of orders $p$, $q$ and $r$, respectively (or any pair of these, see above). By the rigidity of turnover groups in Isom($\bb{H}^{3}$) \cite[Chapter IX.C]{Maskit87}, there is a geodesic plane $\Pi_{T} \subset \bb{H}^{3}$ which is invariant under the action of $T$, and on which $T$ acts as a standard turnover group, with fundamental domain consisting of two $\D(p,q,r)$ triangles. Thus, by the Gauss-Bonnet theorem, the  area  of $\cal{T} := \Pi_{T}/T$ is $2\p(1 - (\frac{1}{p} + \frac{1}{q} + \frac{1}{r}))$, and $\cal{T}$ is a $2$--dimensional space, homeomorphic to a $2$--sphere, with a Riemannian metric of constant curvature $-1$ in the complement of three points and such that  each of these three points has a neighborhood isometric to a hyperbolic cone.  We will refer to such a $\cal{T}$ as a \emph{hyperbolic turnover}, and denote by $f\co \cal{T} \to Q$ the restriction to $\cal{T}$ of the covering map $\bb{H}^{3}/T \to Q$, and we will call $f(\cal{T}) \subset Q$ an \emph{immersion} of a hyperbolic turnover in the hyperbolic $3$--orbifold $Q$. Notice that $\Phi^{-1}(f(\cal{T})) = \bigcup_{\gamma \in \C} \gamma \Pi_{T}$ is a union of geodesic planes in $\bb{H}^{3}$.  We will say in this case that  the immersion is \emph{totally geodesic}. Thus, an immersion of a hyperbolic turnover $f\co \cal{T} \to Q$ is totally geodesic if and only if the map is injective on the level of fundamental groups. We will call the non-orientable $2$--orbifold doubly covered by a turnover a \emph{triangle with mirrored sides}, or often just a \emph{triangle} if the context is clear.  

There are turnovers with spherical and Euclidean structures as well, which are obtained as index two orientation-preserving subgroups of discrete groups generated by reflections in the sides of spherical or Euclidean triangles, respectively. Using the notation of triples to denote the submultiples of $\pi$ in the given triangle, there are four types of spherical turnover: $(2,3,3), (2,3,4), (2,3,5)$, and $(2,2,n)$ for $n \geq 2$.  The first three correspond to the quotient of the $2$--sphere $\bb{S}^{2}$ by the orientation-preserving isometries of the regular tetrahedron, cube, and dodecahedron, respectively. The last is the quotient of $\bb{S}^{2}$ by the dihedral group generated by two order $2$ rotations whose axes meet in an angle of $\pi/n$. If we add to this list both $\bb{S}^{2}$ and the quotient of $\bb{S}^{2}$ by a cyclic group of order $n$, then it is well known that every point in a complete orientable hyperbolic $3$--orbifold has a neighborhood which is isometric to the cone on one of these six types of spherical $2$--orbifolds.

The Euclidean turnovers are $(2,3,6), (2,4,4), (3,3,3),$ and $(2,2,\infty)$. The first three correspond to the doubles of the standard Euclidean triangles, and the last is an open-ended noncompact ``pillowcase'' (i.e., the double of a finite-width infinite half-strip in the plane). The other Euclidean $2$--orbifolds we will consider are the torus and the quotient ``pillow'' obtained from the torus by an involution with four fixed points.

Let $\cal{O}$ be a compact $3$--orbifold, possibly with boundary.  We use the term \emph{Haken ball} to refer to the quotient of a compact $3$--ball by a finite group of isometries.  We say $\cal{O}$ is \emph{irreducible} if every embedded spherical $2$--suborbifold of $\cal{O}$ bounds a Haken ball in $\cal{O}$, and \emph{atoroidal} if every $\pi_{1}$--injective map of a torus, pillow, or Euclidean turnover into $\cal{O}$ is parallel into a boundary component (parallel in the sense that the $2$--orbifold cuts off a product neighborhood of a boundary component). An arbitrary $3$--orbifold is called \emph{atoroidal} if it is diffeomorphic to the interior of an atoroidal $3$--orbifold with boundary. A $2$--orbifold $F$ in $\cal{O}$ is \emph{compressible} if either $F$ is a spherical $2$--orbifold which bounds a Haken ball or if there is a homotopically nontrivial curve in $F$ which bounds a disk quotient (i.e., an \emph{orbifold disk}) in $\cal{O}$. We say $F$ is \emph{incompressible} otherwise. Similarly, there is a relative notion of \emph{$\partial$--compressibility} and \emph{$\partial$--incompressibility} (whose exact definition we will not require).  We call $F$ \emph{essential} if it is incompressible, $\partial$--incompressible, and not parallel into a boundary component. We call a compact irreducible $3$--orbifold \emph{Haken} if it is either a Haken ball, or a turnover crossed with an interval, or if it contains an essential $2$--suborbifold but contains no essential turnover. A compact irreducible $3$--orbifold is called \emph{small} if contains no essential $2$--suborbifolds and has (possibly empty) boundary consisting only of turnovers. The definitions of \emph{essential, Haken,} and \emph{small} extend to arbitrary $3$--orbifolds in the same way as the definition of atoroidal.  Finally, we say a $n$--orbifold is \emph{good} if it is covered by a $n$--manifold. All orbifolds considered in the paper here are assumed to be good.

\section{Analysis of the Complement of the Turnover}\label{S:TurnoverComp}
The idea behind the proof of Theorem \ref{T:ImmersedTurnover} is essentially due to Cooper \cite{Cooper99-1}, that is, to cut the $3$--orbifold along the turnover and then bound the volumes of the complementary pieces. In this way, we will bound from above the volume of the $3$--orbifold.

\subsection{Preliminaries}\label{SS:Prelims}
To begin, let $T \leq \C \cong \p_{1}(Q)$ be a turnover subgroup of the fundamental group of the complete orientable hyperbolic $3$--orbifold $Q$, and let $f\co \cal{T} \to Q$ be the associated isometric immersion.  Let $\Phi\co \bb{H}^{3} \to Q$ be the covering projection. We will require the following well known result. Its proof follows from the proper discontinuity of the action of $\C$ and the Margulis lemma.

\begin{proposition}\label{P:LocallyFinite}
	Let $f\co S \to Q$ be a totally geodesic immersion of any $2$--orbifold $S$ of finite area in a hyperbolic $3$--orbifold $Q$. Let $\C < \Isom(\bb{H}^{3})$ be the fundamental group of $Q$, and let $\C_{S} \leq \C$ be an isomorphic copy of $\p_{1}(S)$ stabilizing a geodesic plane $\Pi \subset \bb{H}^{3}$. Then the collection of planes $\Phi^{-1}(f(S)) = \bigcup_{\gamma \in \C} \gamma\Pi$ is locally finite, in the sense that any compact set $K \subset \bb{H}^{3}$ meets only finitely many planes in the collection.
\end{proposition}

Let $\bigcup_{j \in \cal{J}} S_{j}$ be the disjoint union of the  components of $Q - f(\cal{T})$, and let $\Sigma_{j}$ be a connected component of $\Phi^{-1}(S_{j})$.  Define $\C_{j}$ to be the stablizer of $\Sigma_{j}$ in $\C$, and let $\Lambda(\C_{j})$ denote the limit set of $\C_{j}$, that is, the collection of accumulation points of the orbit space $\C_{j}x$ on the sphere at infinity $S_{\infty}^{2}$ of $\bb{H}^{3}$, where $x \in \bb{H}^{3}$ is any point. Note that $\C_{j} \cong \p_{1}(S_{j})$.  Denote the set $\bb{H}^{3} \cup S_{\infty}^{2}$ by $\overline{\bb{H}}^{3}$. The \emph{convex hull of} $\Lambda(\C_{j})$, which is denoted by $CH (\Lambda(\C_{j}))$, is defined to be the smallest convex set in $\overline{\bb{H}}^{3}$ containing $\Lambda (\C_{j})$ (equivalently, $CH(\Lambda(\C_{j}))$ is the intersection of all the closed half-spaces of $\overline{\bb{H}}^{3}$ which contain $\Lambda(\C_{j})$).  The set $CC(S_{j}) := \Phi(CH (\Lambda(\C_{j})) \cap \bb{H}^{3} )$ is called the \emph{convex core of} $S_{j}$. It is the smallest convex subset of $S_{j}$ which carries its fundamental group.  

Let $C \subset \bb{H}^{3}$ be a closed convex set which is invariant by a Kleinian group $G$, and for this paragraph only denote by $\overline{C}$ the closure of $C$ in $\overline{\bb{H}}^{3}$. Then there is well-defined, $G$--invariant projection $r\co \overline{\bb{H}}^{3} - \overline{C} \to C$ which is given by nearest point retraction. This is defined as follows. If $x \in \overline{\bb{H}}^{3} - \overline{C}$, then take $B(x,\rho)$ to be a ball of radius $\rho$ (or horoball, if $x \in S^{2}_{\infty}$) around $x$ and expand this ball until it meets $C$. This unique first point of intersection is $r(x)$. It is the intersection of a ball (or horoball) around $x$ with a supporting half-plane for the convex set $C$.  The inverse image of the frontier of $C$, $r^{-1}({\mbox{\rm{fr}}}(C)) =  r^{-1}(C) \cap \overline{\bb{H}^{3} - C}$, is just the union, over all $z \in {\mbox{\rm{fr}}}(C)$, of rays beginning at $z$ and which are perpendicular to a supporting half-plane at $z$. Furthermore, if $C$ is not contained in a  proper hyperbolic subspace, then ${\mbox{\rm{fr}}}(C)$ is homeomorphic to $S^{2}_{\infty} - \overline{C}$. The proofs of these facts are all contained in the work of Epstein and Marden \cite[Sections 1.2--1.4]{EpsteinMarden87}. 

In particular, when the convex set in question is $CH(\Lambda(\C_{j}))$, we can use these facts about the nearest point retraction to determine the shape of $\Sigma_{j}$. We will conduct the analysis based on the cardinality of the limit set $\Lambda(\C_{j})$, which contains $0, 1, 2$ or infinitely many points \cite[Chapter II.D]{Maskit87}.  We call a Kleinian group \emph{elementary} or \emph{nonelementary} depending on whether or not it has a finite limit set.  Elementary Kleinian groups are those which contain an abelian subgroup of finite index. 

\subsection{Case: $| \Lambda(\C_{j}) | = 0$}\label{SS:0PtLimSet}
According to the classification of elementary Kleinian groups \cite[Chapter V]{Maskit87}, $\C_{j}$ is a finite group isomorphic to the trivial group, a cyclic group, a dihedral group, or the group of orientation-preserving rigid motions of one of the five Platonic solids. In particular, $\C_{j}$ fixes a point in $z \in \Sigma_{j}$, and $\Sigma_{j}$ is the union of all geodesic rays from $z$ that miss $\Phi^{-1}(f(\cal{T}))$ and all  segments from $z$ to $\Phi^{-1}(f(\cal{T}))$ whose interiors miss $\Phi^{-1}(f(\cal{T}))$. 

The set $\Phi^{-1}(f(\cal{T}))$ is locally finite and its complement in $\bb{H}^{3}$ is homeomorphic to a collection of open balls (the latter fact is most easily seen using the projective ball model of $\bb{H}^{3}$).  Thus, $\overline{\Sigma}_{j}$ is a solid hyperbolic polyhedron $P$ (henceforth, we will use the bar to  denote closure in the ambient space, either $\bb{H}^{3}$ or $Q$). Because $\cal{T}$ has finite diameter and finite area, the closure $\overline{S}_{j}$ of a component $S_{j}$ of $Q - f(\cal{T})$ with finite fundamental group must be compact. Therefore, since $P$ covers $\overline{S}_{j}$ by the action of the finite group $\C_{j}$, $P$ must be a compact polyhedron.  The isoperimetric inequality for hyperbolic $3$--space (e.g. \cite[Section 9]{BurZal}, \cite[Section 6.4]{Chavel93}) implies that the ratio of surface area to volume enclosed is always greater than 2. Applying this to $\Sigma_{j}$ gives $\Vol(\Sigma_{j}) < \Area(\partial \overline{\Sigma}_{j}) / 2$, and this inequality passes to the quotient:

\begin{lemma}\label{L:LimSet0}
	If $S_{j}$ is a component of the complement of the immersed turnover $f(\cal{T})$ in $Q$  and $\p_{1}(S_{j})$ is elementary with empty limit set, then $CC(S_{j})$ is a point and $$\Vol(S_{j}) < {\Area(\partial \overline{S}_{j}) \over 2}.$$
\end{lemma}

\subsection{Case: $| \Lambda(\C_{j}) | = 1$}\label{SS:1PtLimSet}	
In this case, $\C_{j}$ fixes a point on the sphere at infinity $S^{2}_{\infty}$, which we will assume to be the point $\infty$ in the upper half-space model of $\bb{H}^{3}$, and the connected component $\overline{\Sigma}_{j}$ is the union of all geodesics segments from $\infty$ to $\Phi^{-1}(f(\cal{T}))$ \emph{which are contained in} $\overline{\Sigma}_{j}$. It remains to determine what this collection looks like. Again we appeal to the classification of elementary Kleinian groups. In particular, the maximal parabolic subgroup of $\C_{j}$ has rank $1$ or rank $2$. 

\begin{figure}[htbp]
		\begin{center}
		\psfrag{Bounding Plane}{\small{Bounding Plane}}
		\scalebox{1}{\includegraphics{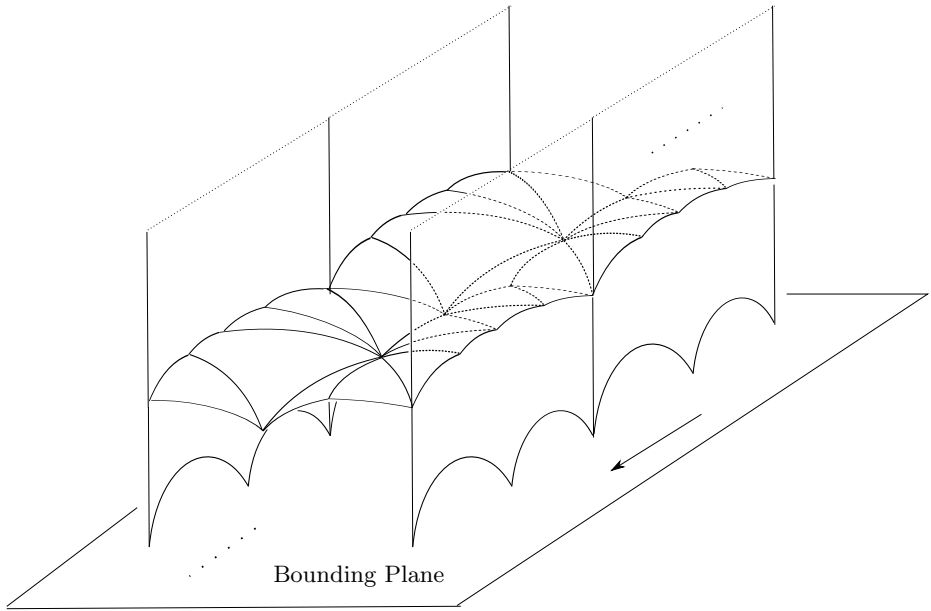}}
		\caption{A vertical ``slab'' cut out of a horoball by $\Phi^{-1}(f(\cal{T}))$, with direction 
			of the action of $\C_{j}$ indicated by the arrow in the bounding plane}
		\label{F:TiledAnnulus1}
	\end{center}
	\end{figure}

We begin with the following observation. In both the rank $1$ and rank $2$ cases, the classification of elementary Kleinian groups provides that $\C_{j}$ stabilizes any horoball centered at $\infty$. We would like to see that $\overline{\Sigma}_{j}$ must project ``vertically'' onto such a horoball. If not, then there must be  vertical planes in $\Phi^{-1}(f(\cal{T}))$ (that is, planes containing the point $\infty$) which cut out some subset of a horoball centered at $\infty$ in such a way that the subset is invariant by $\C_{j}$. For example, $\Sigma_{j}$ might have two vertical ``walls'' which cut out a ``slab'' on which $\C_{j}$ acts by translations. See Figure \ref{F:TiledAnnulus1}, which is meant to illustrate this possibility. Each of the polygonal pieces should be imagined as part of hemisphere that is perpendicular to the bounding plane (i.e., a polygon in a geodesic plane of $\bb{H}^{3}$). 
We observe, however, that this possibility cannot occur for the simple reason (observed above in the previous case) that a compact turnover has finite diameter. This observation, combined with the fact that $\Sigma_{j}$ is connected and stabilized by $\C_{j}$,  proves the following:

\begin{lemma}\label{L:HoroIsotopy}
	If $S_{j}$ is a component of $Q - f(\cal{T})$ and $\p_{1}(S_{j})$ is elementary with limit set a single point, then the closure of a connected component $\Sigma_{j}$ of $\Phi^{-1}(S_{j})$ is isotopic to a horoball centered at the fixed point of $\C_{j} = \Stab_{\C}(\Sigma_{j})$, with the isotopy given by projection toward this fixed point.
\end{lemma}

If the rank of the maximal parabolic subgroup of $\C_{j}$ is $1$, then $\C_{j}$ is isomorphic to either $\bb{Z}$ or a Euclidean $(2,2,\infty)$ turnover group, each of these groups acting on horoballs centered at $\infty$. In either case, we will derive a contradiction.  Observe that a fundamental domain for the action of $\C_{j}$ on a horosphere centered at $\infty$ has infinite area.  At this point, we have not shown that the collection of planes $\Phi^{-1}(f(\cal{T}))$ is connected. Consequently, it could happen that $\partial \overline{\Sigma}_{j}$ does not project onto the whole horosphere.  If the projection \emph{is} onto the whole horosphere, then the isotopy of Lemma \ref{L:HoroIsotopy} maps a fundamental domain for $\C_{j}$ in $\partial \overline{\Sigma}_{j}$ onto a set of infinite area in a fundamental domain for $\C_{j}$ in a horosphere centered at $\infty$. This gives an immediate contradiction, because $f(\cal{T})$ has finite area in $Q$.  But if the mapping is onto a set of \emph{finite} area in a fundamental domain for $\C_{j}$ in a horosphere centered at $\infty$, then we still have a contradiction, as follows.

The subgroup $\C_{j} < \C$ acts on horoballs so that the quotient component $S_{j} = \Sigma_{j} / \C_{j} \subset Q - f(\cal{T})$ is an open convex solid torus or solid pillow. The boundary $\partial \overline{S}_{j}$, which is composed of pieces of $f(\cal{T})$, is a compact set in $Q$ with no boundary. As a consequence, it can have no ends.  We must therefore have that $\partial \overline{S}_{j}$ lifts to a set in a fundamental domain for the action of $\C_{j}$ which has no ends. Now we have that $\Phi^{-1}(f(\cal{T}))$ is a union of geodesic planes which cuts out a collection of open $3$--balls in $\bb{H}^{3}$, and we also have that a fundamental domain for $\C_{j}$ in a horosphere centered at $\infty$ has infinite area. These two facts imply that, if the image (under the above isotopy) of a fundamental domain for $\C_{j}$ in $\partial \overline{\Sigma}_{j}$ has finite area in a fundamental domain in a horosphere, then $\partial \overline{S}_{j}$ must lift to a set with ends in a fundamental domain for $\C_{j}$. But this cannot happen.  We conclude that the rank $1$ case cannot  occur.

\begin{figure}[htbp]
		\begin{center}
		\psfrag{Bounding Plane}{{\small Bounding Plane}}
		\scalebox{1}{\includegraphics{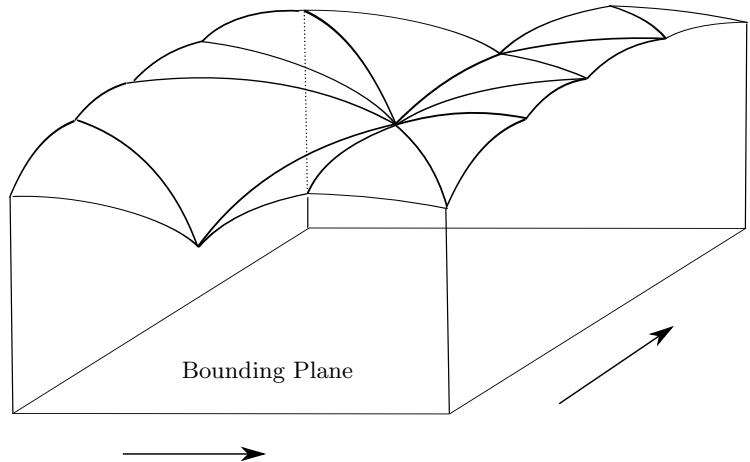}}
		\caption{The fundamental domain for the action of $\C_{j}$ in the rank $2$ case}
		\label{F:TiledCusp1}
		\end{center}
	\end{figure}

If the rank of the maximal parabolic subgroup of $\C_{j}$ is $2$, then by the classification of elementary Kleinian groups, $\C_{j}$ is isomorphic to one of the compact Euclidean turnover groups ($(2,3,6), (2,4,4),$ or $(3,3,3)$), or the fundamental group of a pillowcase or torus. When we project to $Q$, we see that $f(\cal{T})$ cuts out a neighborhood of a finite-volume cusp, which is called \emph{rigid} in the Euclidean turnover case and \emph{non-rigid} otherwise. See Figure \ref{F:TiledCusp1} for a possible illustration in the torus case.

We will eventually see that such a cusp must be rigid (this is, in essence, due to the fact that a hyperbolic turnover has a unique hyperbolic structure). In this case, that is, the compact Euclidean turnover group case, we can bound the volume of a fundamental domain in $\overline{\Sigma}_{j}$ for $\C_{j}$ in terms of the area of (a portion of) its boundary. We let $F$ denote such a fundamental domain, and we choose $F$ to be a union of totally geodesic simplices. The set $\partial F \cap \partial \overline{\Sigma}_{j}$ consists of totally geodesic simplices, each point of which has a unique geodesic segment connecting it to $\infty$ (that is, an observer at $\infty$ can see all of $\partial F \cap \partial \overline{\Sigma}_{j}$). Observe that $\Area(\partial F \cap \partial \overline{\Sigma}_{j}) = \Area(\partial \overline{S}_{j})$. The following result is presumably well known. We provide the argument for completeness.

\begin{lemma}\label{L:GeodesicTriangleVolBound}
 $\Vol(F) < {1 \over 2} \Area(\partial F \cap \partial \overline{\Sigma}_{j})$. 
\end{lemma}

\begin{figure}[htbp]
		\begin{center}
		\psfrag{Bounding Plane}{{\small Bounding Plane}}
		\psfrag{F}{$F$}
		\psfrag{T}{$\triangle$}
		\scalebox{1}{\includegraphics{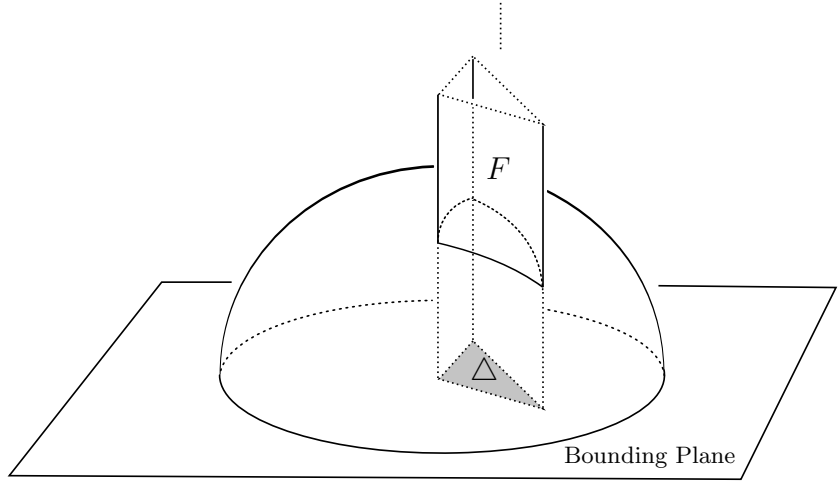}}
		\caption{The region, in upper half-space, above a geodesic triangle}
		\label{F:EuclidVolCalc}
		\end{center}
	\end{figure}

\emph{Proof of \ref{L:GeodesicTriangleVolBound}.}  Since $\partial F \cap \partial \overline{\Sigma}_{j}$ is made up of totally geodesic triangles, it suffices to prove the result when $F$ is just the region above a geodesic triangle in the upper half-space model of $\bb{H}^{3}$.  See Figure \ref{F:EuclidVolCalc}.  After applying an isometry of $\bb{H}^{3}$, we may assume that the hemisphere representing the geodesic plane containing our triangle has Euclidean radius 1. Let $\triangle$ be the triangle in the bounding plane which is the Euclidean orthogonal projection of  $\partial F \cap \partial \overline{\Sigma}_{j}$. Note that the equatorial disk of our unit hemisphere represents a copy of the projective model for $\bb{H}^{2}$, and $\triangle$ is an isometric copy of $\partial F \cap \partial \overline{\Sigma}_{j}$ (in the projective model) under the projection.  We have the following calculation
\[	\Vol(F) = \iint\limits_{(x,y) \in \triangle} \left( \int_{\sqrt{1 - (x^{2}+y^{2})}}^{\infty} {dz \over z^{3}} \right)
	 \, dx \, dy 
\]
\[	= {1 \over 2} \iint\limits_{(x,y) \in \triangle} {dx \, dy \over 1 - (x^{2}+y^{2})}
	< {1 \over 2} \iint\limits_{(x,y) \in \triangle} {dx \, dy \over (1 - (x^{2}+y^{2}))^{{3 \over 2}}}
	= {1 \over 2} \Area(\partial F \cap \partial \overline{\Sigma}_{j}).
\]\\
The final equality is obtained by recognizing the integrand of the last double integral as the area element for the projective model of $\bb{H}^{2}$. This proves the lemma. \hfill  \fbox{\ref{L:GeodesicTriangleVolBound}}\\

Observe that this result applies to any finite-volume cusp neighborhood which is cut out from $Q$ by $f(\cal{T})$. By Corollary \ref{C:NoNonTOs}, however, the non-rigid cusp neighborhoods will not occur. In particular, we have the following:

\begin{lemma}\label{L:LimSet1Rigid}
	If $S_{j}$ is a component of $Q - f(\cal{T})$ and $\p_{1}(S_{j})$ is elementary with $1$ limit point, then $S_{j}$ is a finite-volume, rigid cusp neighborhood in $Q$, and $$\Vol(S_{j}) < { \Area(\partial \overline{S}_{j}) \over 2}.$$
\end{lemma}

\subsection{Case: $| \Lambda(\C_{j}) | = \infty$}\label{SS:InfLimSet}
Since $\Lambda(\C_{j})$ is contained in one ``side'' of $S^{2}_{\infty}$ determined by any plane $\Pi \subset \Phi^{-1}(f(\cal{T}))$ which contains a facet of $\overline{\Sigma}_{j}$, it follows that $CH(\Lambda(\C_{j})) \neq \bb{H}^{3}$. So $CH(\Lambda(\C_{j}))$ projects to $CC(S_{j})$ in $Q$ as either a totally geodesic hyperbolic $2$--suborbifold or a hyperbolic $3$--suborbifold with a hyperbolic $2$--orbifold in its boundary which separates $f(\cal{T})$ from $CC(S_{j}$) \cite[Proposition 8.5.1]{Thurstonnotes}, \cite[Theorem 1.12.1]{EpsteinMarden87}. The flow induced by the gradient of the function which gives the distance from a point to $CH(\Lambda(\C_{j}))$ determines a product structure (e.g. \cite[Appendix A]{Agol04-1})
	$$\overline{S}_{j} - CC(S_{j}) \cong \partial \overline{S}_{j} \times [0,1).$$
We make some remarks to give an idea of what is going on here. Suppose that $CH(\Lambda(\C_{j}))$ is a geodesic plane.  If $\C_{j}$ contains no elliptic element of order two whose axis is contained in this plane, i.e., if $\C_{j}$ is an orientation-preserving Fuchsian group, then we obtain the product structure
	$$\overline{S}_{j} \cong \partial \overline{S}_{j} \times [0,1],$$
with $CH(\Lambda(\C_{j})) = \partial \overline{S}_{j} \times \{1/2\}$. In the case that $\C_{j}$ \emph{does} contain an elliptic element of order two whose axis is contained in $CH(\Lambda(\C_{j}))$, then $\Sigma_{j}$ is a neighborhood of a geodesic plane, and this order two element of $\C_{j} = \Stab_{\C}(\Sigma_{j})$ exchanges the two halves of this neighborhood. Such an order two element  makes $CC(S_{j})$ non-orientable as a $2$--orbifold, and $1$--sided as a $2$--\emph{sub}orbifold. In particular, a regular neighborhood of $CC(S_{j})$ in this case will have boundary a $2$--sided $2$--suborbifold $F \subset Q$ which doubly covers  $CC(S_{j})$, and we have a homeomorphism $\partial \overline{S}_{j} \times \{t\} \cong F$ for any $t \in [0,1)$.  Finally, when $CC(S_{j})$ is $3$--dimensional and $F$ is the hyperbolic $2$--orbifold in $\partial CC(S_{j})$ which separates the convex core of $S_{j}$ from $f(\cal{T})$, we have the homeomorphism $\partial \overline{S}_{j} \times \{t\} \cong F$ for any $t \in [0,1)$.

Consider the projection $r_{\partial}\co \partial \overline{\Sigma}_{j} \to CH(\Lambda(\C_{j}))$ given by restricting $r$. Because $r$ is distance decreasing \cite[Lemma 1.3.4]{EpsteinMarden87}, $r_{\partial}$ is area decreasing. So the image $\Phi \circ r_{\partial}(\partial \overline{\Sigma}_{j})$ in $Q$ is a collection of hyperbolic $2$--orbifolds with total area less than twice the area of $\cal{T}$. Define $R_{j}$ to be the image under the projection $\Phi$ of the product region between $\partial \overline{\Sigma}_{j}$ and $r_{\partial}(\partial \overline{\Sigma}_{j})$, the possibilities for which were described in the previous paragraph.  So we have 
	$$R_{j} = S_{j} - CC(S_{j}).$$ 
(For reasons which will be obvious in Section \ref{S:Isop}, this product region will be called a \emph{room}). In Section \ref{S:NbhdHyp}, we will see that it is impossible for $\Phi \circ r_{\partial}(\partial \overline{\Sigma}_{j})$ to be anything other than a collection of hyperbolic turnovers or mirrored hyperbolic triangles (Corollary \ref{C:NoNonTOs}). Thus $R_{j}$ is homeomorphically a collection of  turnovers crossed with an interval. We will also see, in Section \ref{S:Isop}, that in this case we have the volume bound
	$$\Vol(R_{j}) < {H \over 2} \Area(\partial \overline{S}_{j}),$$
where $H=1.199678...$ is the positive solution of $x=\coth x$. We summarize these remarks as follows:

\begin{lemma}\label{L:LimSetInf}
	Suppose $S_{j}$ is a component of $Q - f(\cal{T})$ and $\p_{1}(S_{j})$ is non-elementary. Then the region $R_{j}$ between $\partial \overline{S}_{j}$ and $CC(S_{j})$ is a product of a turnover with an open interval, and $\Vol(R_{j}) <  {H \over 2} \Area(\partial \overline{S}_{j})$, for $H$ described above.
\end{lemma}

\subsection{Case: $| \Lambda(\C_{j}) | = 2$}\label{SS:2PtLimSet}
\begin{figure}[htbp]
		\begin{center}
		
		\psfrag{Sig_j}{{\small $\Sigma_{j}$}}
		\psfrag{Lam(Gam_j)}{{\small $\Lambda(\C_{j})$}}
		\psfrag{projects to}{{\small projects to}}
		\psfrag{core curve of}{{\small core curve of}}
		\psfrag{solid torus/pillow}{{\small solid torus/pillow}}
		\scalebox{1}{\includegraphics{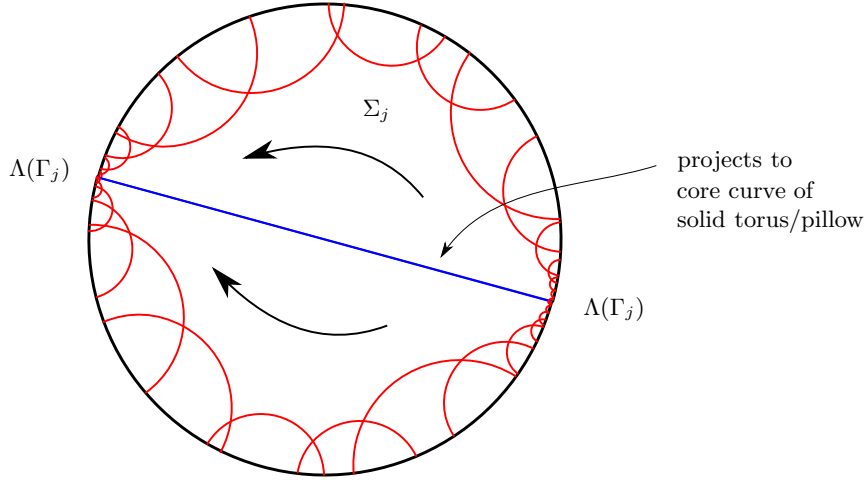}}
		\caption{Schematic diagram for the two-point limit set case}
		\label{F:2ptlimitset}
		\end{center}
	\end{figure}

Let $\Lambda(\C_{j}) = \{x,y\} \subset S^{2}_{\infty}$. Then the convex hull is the geodesic $\tilde{l}$ connecting $x$ and $y$, and the support planes for this convex hull consist of all planes containing $\tilde{l}$. It follows that $\Sigma_{j}$ is a neighborhood of $\tilde{l}$, and $\C_{j}$ stabilizes this neighborhood, acting by translations or as an infinite dihedral group along $\tilde{l}$, and possibly also by rotations around $\tilde{l}$. Figure \ref{F:2ptlimitset} provides a schematic/lower-dimensional illustration of  a collection of $\C_{j}$--invariant planes in $\bb{H}^{3}$ cutting out a neighborhood of a stabilized geodesic.

Thus $f(\cal{T})$ cuts out a neighborhood of the orbifold geodesic $l= \Phi(\tilde{l})$ in $Q$.  In this circumstance, we can put a complete hyperbolic structure on $Q_{l} = Q - l$ \cite[Theorem 1.2.1]{Kojima98-1}, and so we are left with a new hyperbolic $3$--orbifold, with a new cusp $C$, which again contains a new hyperbolic turnover $g(\cal{T})$.  This new cusp lies in a component of $Q_{l} - g(\cal{T})$ whose fundamental group has limit set consisting either of $1$ point or infinitely many points. This cusp admits deformations (for example, simply fill the cusp back in to obtain $Q$), so in the former case we would have that $g(\cal{T})$ cuts out a non-rigid cusp neighborhood in $Q_{l}$. We have remarked in Section \ref{SS:1PtLimSet} that such a phenomenon does not occur for an immersed turnover (see Lemma \ref{L:LimSet1Rigid}), and so we must be in the latter case. Furthermore, we have remarked that if the fundamental group of this component has infinite limit set, then Corollary \ref{C:NoNonTOs} tells us its convex core must be a $3$--dimensional orbifold with boundary a collection of embedded totally geodesic hyperbolic turnovers $\cal{S}$ separating it from the rest of $Q_{l}$. 

We observe that any such turnover boundary component will remain $\pi_{1}$--injective under the  hyperbolic Dehn filling on the cusp $C$ which yields our original orbifold $Q$. This is by the Equivariant Loop Theorem \cite{MeeksYau80-1}, which implies that any embedded $2$--orbifold in a good $3$--orbifold which is not $\pi_{1}$--injective has a compressing orbifold disk. If some elements of $\pi_{1}(\cal{S})$ were made trivial by filling $C$ to retrieve $Q$, then there would have to be compressing orbifold disks to achieve this for the image of $\cal{S}$ in $Q$. But a turnover admits no compressing orbifold disks, because it admits no essential simple closed curves.  We therefore have a non-empty collection of embedded totally geodesic turnovers in $Q_{l}$ which separate the cusp $C$ from the immersed turnover $g(\cal{T})$, and this collection survives any hyperbolic Dehn surgery performed on $C$.  In particular, the collection survives if we fill $C$ to retrieve $Q$. But now the collection $\cal{S}$ geometrically isolates $C$ from $g(\cal{T})$.  This implies that the geometry of the immersion $g\co \cal{T} \to Q_{l}$ does not change under fillings of $C$.  In particular, there is a hyperbolic $3$--orbifold $V$ which has some totally geodesic boundary components, which contains a totally geodesic immersion $h(\cal{T})$ of the turnover $\cal{T}$, and for which the pair $(V, h(\cal{T}))$ isometrically embeds in both $(Q, f(\cal{T}))$ and $(Q_{l}, g(\cal{T}))$.  But this is a contradiction, because we began by assuming our immersed turnover $f\co \cal{T} \to Q$ cut out tubular neighborhood of an orbifold geodesic, which cannot contain any totally geodesic hyperbolic turnovers.  This contradiction rules out the two-point limit set case:

\begin{lemma}\label{L:LimSet2}
	If $S_{j}$ is a component of $Q - f(\cal{T})$, then $\p_{1}(S_{j})$ can not be elementary with a two-point limit set.
\end{lemma}

\section{The ``Geodesic Floor'' Isoperimetric Inequality}\label{S:Isop}
In order to prove Lemma \ref{L:LimSetInf}, we need a tool for bounding the volume between $S_{j}$ and $CC(S_{j})$, when $\p_{1}(S_{j})$ is non-elementary and either $CC(S_{j})$ is $2$--dimensional and totally geodesic or $CC(S_{j})$ is $3$--dimensional and $\partial CC(S_{j})$ is totally geodesic. This will be provided by the next theorem. First, we need some definitions and notation.

Let $F$ be a measurable subset of finite area in a geodesic plane of $\bb{H}^3$.  Choose one closed half-space $H = \bb{H}_{+}^{3} \subset \bb{H}^{3}$ of the plane containing $F$. For any $z \in F$, let $l_{z}$ denote the geodesic ray perpendicular to $F$ that begins at $z$ and points into $H$. A \emph{graph} over $F$ is the image of a function $g\co F \to H$ which maps $z \in F$ to the point on $l_{z}$ at (hyperbolic) distance $\varphi (z)$ from $z$, where $\varphi\co F \to \bb{R}_{\geq 0}$ is a non-negative function.  Let $C$ be a graph over $F$ and $R := R(F,C)$ the region trapped between $F$ and $C$. We will call $F$, $C$ and $R$ the \emph{floor}, the \emph{ceiling} and the \emph{room}, respectively. 

 Let $S(r)$ denote the graph of constant height $r$ over $F$, and let $B(r)$ denote the region trapped between $S(r)$ and $F$. We will refer to such an $S$ (respectively, $B$) as a \emph{nice ceiling} (respectively, \emph{nice room}).

\begin{theorem}\label{T:IsopIneq}
(Geodesic Floor Isoperimetric Inequality) Let $C$ be the ceiling for the graph $g\co F \to H$ of a function $\varphi \co F \to \bb{R}_{\geq 0}$ which, for simplicity, is assumed to be differentiable almost everywhere. Let $R(F,C)$ be the associated room. Let $B$ be the nice room over $F$ with $\Vol(B) = \Vol(R)$, and let $S$ be the ceiling of $B$. Then $\Area(C) \geq \Area(S)$.
\end{theorem}

\emph{Proof of \ref{T:IsopIneq}.} 
We take the metric on $\bb{H}^{3}$ given by $dh^{2} + \cosh^{2}h(dr^{2} + \sinh^{2} r d\th^{2})$, where $r$ and $\th$ describe polar coordinates in the plane containing $F$ and $h$ is the positive hyperbolic distance into $H$ from $F$. We also denote the area form $\sinh r \, dr d\th$ on the plane containing $F$ by $dA$. Let $g\co F \to H$ be given, and $C$ the image of $g$. The volume of the room $R$ is given by
	\begin{align}\label{E:Vol(R)}
		V := \Vol(R) &= \iint\limits_{F} \int_{0}^{g(r,\th)} \cosh^{2}h \, dh \, dA\\ 
		&= \iint\limits_{F} {1 \over 4} \left( \sinh 2g(r,\th) + 2g(r,\th) \right) \, dA,\notag
	\end{align}
and since the volume of the nice room $B$ is the same, we have the constant $H$ defined implicitly by
	\begin{equation}\label{E:Vol(B)}
		\Vol(B) = \iint\limits_{F} {1 \over 4} \left( \sinh 2H + 2H \right) \, dA = V.
	\end{equation}
A calculation yields the following formula and lower bound for the area of $C$
	\begin{align}\label{E:Area(C)}
		\Area(C) &= \iint\limits_{F} \cosh g(r,\th) \sqrt{\left( g_{r}^{2}(r,\th) + 
			\cosh^{2} g(r,\th) \right) \sinh^{2}r + g_{\th}^{2}(r,\th)} \, dr d\th\\ \notag
			&\geq \iint\limits_{F} \cosh g(r,\th) \sqrt{\cosh^{2}g(r,\th) \sinh^{2}r} \, drd\th\\ \notag
			&= \iint\limits_{F} \cosh^{2} g(r,\th) \, dA. \notag
	\end{align}
The constant height $H$ ceiling $S$ of the nice room $B$ has area
	\begin{equation}\label{E:Area(S)}
		\Area(S) = A_{F} \cosh^{2}H = A_{F}(\sinh^{2}H + 1),
	\end{equation}
where $A_{F} = \Area(F)$. Because $H$ is a constant, (\ref{E:Vol(B)}) can be rewritten as 
	\begin{equation}\label{E:VintermsofH}
		\sinh2H + 2H = 4V / A_{F}. 
	\end{equation}
We therefore have
	\begin{equation}
 		{{4V \over A_{F}} - 2H \over 2\cosh H} = {\sinh 2H \over 2 \cosh H} = {2 \sinh H \cosh H 
		\over 2 \cosh H} = \sinh H, \notag
	\end{equation}
which implies
	\begin{equation}
		A_{F}(\sinh^{2} H + 1) = {A_{F}^{2}\cosh^{2} H + (2V - H A_{F})^{2} \over A_{F} \cosh^{2} H}. 
		\notag
	\end{equation}
Substituting (\ref{E:Area(S)}) into the above yields
	\begin{equation}
		\Area(S) = {A_{F} \Area(S) + (2V - H A_{F})^{2} \over \Area(S)}, \notag
	\end{equation}
or, by the quadratic formula,
	\begin{equation}\label{E:Amin}
		\Area(S) = {A_{F} + \sqrt{A_{F}^{2} + 4(2V - H A_{F})^2} \over 2}.
	\end{equation}
We would like to show that the right-hand side of (\ref{E:Amin}) is always less than or equal to $\iint\limits_{F} \cosh^{2} g \, dA$, and therefore, by (\ref{E:Area(C)}), always less than or equal to $\Area(C)$. We begin by observing that 
	\begin{equation}\label{E:H>Avg(g)}
		\iint\limits_{F} H \, dA \geq \iint\limits_{F} g \, dA,
	\end{equation}
for otherwise equations (\ref{E:Vol(R)}) and (\ref{E:Vol(B)}) would imply
	\begin{equation}
		{1 \over A_{F}} \iint\limits_{F} \sinh 2g \, dA < \sinh 2H < 
		\sinh 2 \left( {1 \over A_{F}} \iint\limits_{F} g \, dA \right), \notag
	\end{equation}
which would violate Jensen's inequality, because $x \mapsto \sinh 2x$ is a convex function on $[0,\infty)$. Multiplying both sides of (\ref{E:H>Avg(g)}) by $2$ and adding $\iint\limits_{F} \sinh 2g \, dA$ to both sides, we can apply (\ref{E:Vol(R)}) to obtain
	\begin{equation}
		\iint\limits_{F} 2\sinh g \cosh g \, dA = \iint\limits_{F} \sinh 2g \, dA \geq 4V - 2HA_{F}. \notag
	\end{equation}
Now applying the Cauchy--Schwarz inequality to the left-most integral, we obtain
	\begin{equation}
		\sqrt{\iint\limits_{F} 2 \sinh^{2} g \, dA} \sqrt{\iint\limits_{F} 2 \cosh^{2} g \, dA} \geq 
		4V - 2HA_{F}. \notag
	\end{equation}
The right-hand side of the above inequality is equal to $A_{F}\sinh2H$, which is non-negative. We may therefore square both sides and apply a hyperbolic trigonometric identity to obtain
	\begin{equation}
		\left( \iint\limits_{F} \cosh 2g - 1\, dA \right) \left( \iint\limits_{F} \cosh 2g + 1 \, dA \right) 
		\geq 4(2V - H A_{F})^{2}, \notag
	\end{equation}
and expanding the left-hand side of the above gives
	\begin{equation}
		\left( \iint\limits_{F} \cosh 2g \, dA \right)^{2} - \left( \iint\limits_{F} \, dA \right)^{2} 
		\geq 4(2V - H A_{F})^{2}. \notag
	\end{equation}
We rewrite this as 
	\begin{equation}
		\left( \iint\limits_{F} \cosh 2g + 1 - 1 \, dA \right)^{2} 
		\geq A_{F}^{2} + 4(2V - H A_{F})^{2}, \notag
	\end{equation}
which is equivalent to
	\begin{equation}
		\left( 2 \iint\limits_{F} \cosh^{2} g \, dA  - A_{F} \right)^{2} 
		\geq  A_{F}^{2} + 4(2V - H A_{F})^{2}. \notag
	\end{equation}
Taking square roots, adding $A_{F}$ to both sides, dividing by $2$, and applying (\ref{E:Area(C)}) and (\ref{E:Amin}) to the result finally yields
	\begin{align}
		\Area(C) &\geq \iint\limits_{F} \cosh^{2} g \, dA \\ \notag 
		&\geq {A_{F} + \sqrt{A_{F}^{2} + 4(2V - H A_{F})^2} \over 2} \\ \notag
		&= \Area(S), \notag
	\end{align}
which is what we wanted to show.  \hfill  \fbox{\ref{T:IsopIneq}}\\

In Section \ref{S:Calcs}, we will see that if $B$ is a nice room (over any floor) and $S$ is its ceiling, then 
	$${\Area(S) \over \Vol (B)} > {2 \over H},$$
where $H=1.199678...$ is the positive solution of $x = \coth x$. As a result, we have the following corollary. It implies Lemma \ref{L:LimSetInf}.

\begin{corollary}\label{C:RoomVolBound}
	Let $R(F,C)$ be as in Theorem \ref{T:IsopIneq}. Then $\Vol(R) < {H \over 2} \Area(C)$. 
\end{corollary}

\section{A Hyperbolic ``Neighborhood'' of an Immersed Turnover}\label{S:NbhdHyp}

The goal of this section is to prove the first part of item $(5)$ from Theorem \ref{T:ImmersedTurnover}, i.e., that an immersed turnover $f\co \cal{T} \to Q$ in a complete orientable hyperbolic $3$--orbifold is contained in a small hyperbolic $3$--suborbifold with (possibly empty) totally geodesic boundary.  Recall that a compact $3$--orbifold is \emph{small} if it is irreducible, contains no essential orientable $2$--suborbifold and has (possibly empty) boundary consisting only of turnovers. The uniformization theorem for small $3$--orbifolds \cite{BLP05-1} says that a small $3$--orbifold has a geometric structure.  We continue to denote the covering projection by $\Phi\co \bb{H}^{3} \to Q$, a connected component of $\bb{H}^{3} - \Phi^{-1}(f(\cal{T}))$ by $\Sigma_{j}$, and the stabilizer of such a component by $\C_{j}$.
  
Let $\cal{N}(W)$ denote a closed regular neighborhood of a subset $W \subset Q$. Consider the union of $\cal{N}(f(\cal{T}))$ and all of the components of $Q - f(\cal{T})$ whose fundamental group has finite limit set. From this union, remove small embedded open cusp neighborhoods. Call the result of this construction $N$. It is a compact $3$--orbifold, and $\partial N$ is non-empty if and only if there is a component of $Q - f(\cal{T})$ with infinite fundamental group that is not a neighborhood of an orbifold geodesic. 

\begin{theorem}\label{T:RegNbhdIsHyp}
$N$ is a small $3$--orbifold. Moreover, the geometric structure on $N$ is hyperbolic.
\end{theorem}

\emph{Proof of \ref{T:RegNbhdIsHyp}.}  We begin by showing that $N$ is irreducible. Consider the covering space $\tilde{N}$ of $N$ consisting of a connected component of $\Phi^{-1}(N)$. It is a noncompact manifold with boundary. We will show that $\tilde{N}$ is irreducible and therefore, since $\tilde{N}$ covers $N$, that $N$ is irreducible as well \cite[Remark following Theorem 3.23]{BMP03-1}.

Let $Y$ be an embedded sphere in $\tilde{N}$. We will show that $Y$ bounds a ball. To do this, it is helpful to understand what $\tilde{N}$ looks like. Of course, we have that one component $X$ of $\Phi^{-1}(f(\cal{T}))$ will be contained in $\tilde{N}$.   The complement $\bb{H}^{3} - \Phi^{-1}(f(\cal{T}))$ consists of precompact open balls and non-precompact open balls (this is best visualized using the projective ball model for $\bb{H}^{3}$). The stabilizer for any one of these precompact open balls must be finite because the covering action on $\bb{H}^{3}$ is discrete, and therefore each such ball covers a component of $Q - f(\cal{T})$ with finite fundamental group.  In particular, each open ball in $\bb{H}^{3} - X$ with finite volume is added to $X$ as part of the  construction of $\tilde{N}$. 

We remark that we have not yet ruled out the possibility that $f(\cal{T})$ cuts out a solid torus or pillow from $Q$ (Lemma \ref{L:LimSet2}), for this result depends on both Lemma \ref{L:LimSetInf} and the theorem we are in the middle of proving.  As a consequence, we must include this possibility in the present analysis. We therefore also must add to $X$ (as part of the construction of $\tilde{N}$) any non-precompact ball whose closure meets $X$ and which covers a component of $Q - f(\cal{T})$ that looks like a solid torus or pillow.  These non-precompact balls just look like infinite solid cylinders. 

If a component of $\bb{H}^{3} - X$ whose closure meets $X$ is isotopic to a horoball centered at some $z \in S^{2}_{\infty}$, then we also add this component minus a smaller horoball around $z$  to $X$ as part of $\tilde{N}$.  These pieces of $\tilde{N}$ look like ``horo-slabs,'' and are the truncated developments in $\bb{H}^{3}$ of components of $Q - f(\cal{T})$ that look like cusp neighborhoods.  

To finally complete the construction of $\tilde{N}$, we consider a component of $Q - f(\cal{T})$ whose fundamental group has infinite limit set $\Lambda$. Recall, from the discussion in the first paragraph of Section \ref{SS:InfLimSet}, that there is a product region between $X$ and $\partial CH(\Lambda)$ (or between $X$ and $CH(\Lambda)$, if the convex hull is $2$--dimensional). We add on the ``half'' of this product region which meets $X$ to complete the construction of $\tilde{N}$. It follows from this construction that $\partial \tilde{N}$ consists solely of embedded noncompact surfaces in $\bb{H}^{3}$ (corresponding either to horospheres or to the ``halfway mark'' for the flow onto $CH(\Lambda(\C_{j}))$ for  a component $\Sigma_{j}$ of $\bb{H}^{3} - \Phi^{-1}(f(\cal{T}))$ whose stabilizer has infinite limit set).  
  
Returning to our embedded sphere $Y$, we observe that $Y$ bounds a $3$--ball $B$ in $\bb{H}^{3}$.  If $B$ is contained in $\tilde{N}$, then we are done. If $B \not\subset \tilde{N}$, then some interior point $z \in B$ is not contained in $\tilde{N}$. There are but two possibilities for the location of $z$: it is contained in a horoball region of $\bb{H}^{3} - \tilde{N}$ or it is contained in one of the regions $\Sigma_{j}$ of $\bb{H}^{3} - \Phi^{-1}(f(\cal{T}))$ whose stabilizer has infinite limit set.  A properly embedded noncompact surface separates $\tilde{N}$ from any such region. Since $Y \subset \tilde{N}$, we therefore have that a properly embedded noncompact surface is contained in $B$. This contradiction implies that we must have $B \subset \tilde{N}$, and we conclude that $N$ is irreducible.  

We must show that $N$ contains no essential $2$--orbifolds. We first prove two lemmata.

\begin{lemma}\label{L:NAtoroidal}
$N$ is atoroidal. 
\end{lemma}

\emph{Proof of \ref{L:NAtoroidal}.}  We have just seen that $N$ is irreducible, so in order to show that $N$ is atoroidal it is enough to show that an embedded torus, pillowcase, or Euclidean turnover in $N$ which is incompressible is also $\partial$--parallel. So assume that $S \subset N$ is such a $2$--suborbifold, and suppose $S$ is not parallel to a boundary component of $N$. Observe that as a $2$--suborbifold of the orbifold $Q \supset N$, $S$ must either compress or be parallel into a cusp of $Q$, because $Q$ is atoroidal.  

In the latter case, the product region determined by the isotopy of $S$ into a cusp neighborhood of $Q$ exhibits the isotopy of $S$ into the boundary of $N$, because this cusp neighborhood must have been added to $Q - f(\cal{T})$ in the construction of $N$.  The fact that this isotopy remains in $N$ follows as in the irreducibility argument above, for the lift of the isotopy to the universal cover is a product region between two horospheres centered at the same point on $S_{\infty}^{2}$.  If the isotopy of $S$ does not remain in $N$, then there is be a properly embedded noncompact surface (which is not a horosphere) in this product region. This is impossible. 

If, on the other hand, $S$ compresses in $Q$, then $S$ must be a torus or pillowcase ($S$ cannot be a Euclidean turnover because an embedded circle on a turnover always bounds an orbifold disk in the turnover).  In this case $S$ either bounds a solid orbifold torus or solid orbifold pillow, or $S$ is contained in a Haken ball in $Q$ (e.g. \cite[Proposition 3.14]{BMP03-1}).  If $S$ is contained in a Haken ball in $Q$, then it is homotopically trivial, and again consideration of the construction of $N$ allows us to conclude that $S$ compresses in $N$. If $S$ bounds a solid orbifold torus/pillow $V$, and if $V$ is not contained in $N$, then $V$ must contain a component of $Q - N$. Again, as in the irreducibility argument, we have a contradiction.  This proves the lemma.    \hfill  \fbox{\ref{L:NAtoroidal}}\\

\begin{lemma}\label{L:HypTOBndry||}
If $g\co S \to N$ is an embedded hyperbolic turnover, then $S$ is $\partial$--parallel.
\end{lemma}
\emph{Proof of \ref{L:HypTOBndry||}.}  In the orbifold $Q$, $S$ is isotopic to a totally geodesic embedding or to the double cover of mirrored triangle, which we may as well call $S$ in either case. Consider any curve $\cal{C} \subset S \cap f(\cal{T})$. Because $\cal{C}$ is in the intersection of two totally geodesic orbifolds, it is a geodesic. As a result, $\cal{C}$ cannot be simple, because a geodesic on a hyperbolic turnover must intersect itself. However, we observe that $f^{-1}(\cal{C})$ must be embedded in $\cal{T}$, because otherwise we would have an element of $\p_{1}(\cal{T})$, corresponding to the curve $f^{-1}(\cal{C})$,  taking one lift $\tilde{S}_{1}$ of $S$ to another lift $\tilde{S}_{2}$ with $\tilde{S}_{1} \neq \tilde{S}_{2}$ but $\tilde{S}_{1} \cap \tilde{S}_{2} \neq \emptyset$.  Because $S$ is embedded, this is not possible and we conclude that $f^{-1}(\cal{C})$ is embedded in $\cal{T}$. But an embedded curve on a turnover bounds an orbifold disk, and this implies that the geodesic $\cal{C}$ is contractible, which is a contradiction.  This contradiction also holds in $N$ because $S$ was made totally geodesic in $Q$ by an isotopy. In particular, $S$ is disjoint from $f(\cal{T})$ in $N$, and since we know what the complementary pieces of $N - f(\cal{T})$ look like, the conclusion follows.  \hfill \fbox{\ref{L:HypTOBndry||}}\\

We can use these lemmata to show that $N$ has a hyperbolic structure. To do this, we use the turnover splitting of a $3$--orbifold \cite[Theorem 4.8]{BMP03-1}. This says that a compact, irreducible, atoroidal $3$--orbifold $\cal{O}$ contains a maximal (possibly empty) collection $\frak{h}$ of essential, pairwise disjoint, non-parallel, hyperbolic turnovers which is unique up to isotopy and such that every component of $\cal{O}$ split along $\frak{h}$ is either small or Haken.  Lemma \ref{L:NAtoroidal} shows that $N$ is atoroidal, and Lemma \ref{L:HypTOBndry||} implies that the turnover splitting collection of $N$ is empty. So $N$ is either small or Haken. As we claim, it will turn out that $N$ is small. 

Now suppose that $g\co S \to N$ is a proper embedding of a $\partial$--incompressible, incompressible, orientable $2$--orbifold in $N$. Make $g(S)$ transverse to $f(\cal{T})$ in $N$, and consider any curve $\cal{C} \subset f(\cal{T}) \cap g(S)$. As in the hyperbolic turnover case, we must have that $f^{-1}(\cal{C})$ is simple in $\cal{T}$, for otherwise a multiple point would give rise to a contradiction that $g$ is an embedding. In particular, $f^{-1}(\cal{C})$ bounds an orbifold disk $D_{\cal{T}}$ in $\cal{T}$, and therefore $\cal{C}$ is contractible along $f(D_{\cal{T}})$ in $N$. Now since $g(S)$ is essential, the map $g$ is $\pi_{1}$--injective, and so $\cal{C}$ bounds a (possibly non-embedded) disk $D_{S}$ in $g(S)$. By choosing $\cal{C}$ which has innermost preimage in $\cal{T}$, we can use the irreducibility of $N$ to obtain a map of an orbifold ball $B$ into $N$ with boundary $f(D_{\cal{T}}) \cup D_{S}$. We may then use $B$ to remove $\cal{C}$ from $f(\cal{T}) \cap g(S)$ and thus decrease the number of curves of intersection of the two orbifolds. After a finite number of steps, we see that $g(S)$ can be made disjoint from $f(\cal{T})$. But we know what the components of $N - f(\cal{T})$ look like. Namely, these components are either $\partial$--parallel product regions, orbifold balls, solid tori, or solid pillows. We conclude that $g(S)$ is $\partial$--parallel. This proves that $N$ contains no essential, orientable $2$--suborbifolds.

As we observed above, $N$ is either small or Haken. Because $N$ is irreducible and contains no embedded, orientable, essential $2$--suborbifold, consideration of the construction of $N$ shows that it is small \emph{exactly} when there is no component $\cal{Q} \subset Q - f(\cal{T})$ satisfying any of the following conditions:
\begin{enumerate}
	\item $\cal{Q}$ is a non-rigid cusp neighborhood, or
	\item $CC(\cal{Q})$ is a non-rigid hyperbolic $2$--orbifold, or
	\item $\partial CC(\cal{Q})$ contains a non-rigid hyperbolic $2$--orbifold which separates 
		$CC(\cal{Q})$  from $f(\cal{T})$.
\end{enumerate}
\smallskip
Any such $\cal{Q}$ will produce a non-turnover boundary component in the construction of $N$, and thus prevent $N$ from being small. So to show that $N$ is small, we suppose for the sake of contradiction that there is such a $\cal{Q}$. In this case, we see that $N$ satisfies the following theorem of Dunbar \cite{Dunbar88-1} (note that $N$ contains no bad $2$--suborbifold, i.e., $N$ is ``abad,'' because $Q \supset N$ can contain no such $2$--suborbifold):

\begin{theorem}\label{T:DunbarHeirarchy}
Let $\cal{O}$ be a smooth, compact, connected, irreducible, abad, orientable $3$--orbifold in which every non-spherical turnover is boundary-parallel. If $\partial \cal{O}$ has a component which is not a turnover, then $\cal{O}$ has a ``strong hierarchy,'' i.e., $\cal{O}$ can be decomposed into orbifold balls and thick turnovers by repeated cutting along $2$--sided, essential $2$--suborbifolds.
\end{theorem}

But we have already seen that $N$ contains no $2$--suborbifolds for such a hierarchy. Since $N$ is neither an orbifold ball nor a thick turnover (the latter because $\cal{T}$ is not embedded), we conclude that there can be no such component $\cal{Q}$ as above. Therefore, $N$ is small, and either $\mbox{\rm{int}}(N) \cong Q$ or $N$ is obtained from a component $Q'$ of $Q$ cut along a collection of totally geodesic turnovers and $1$--sided triangles with mirrored sides and possibly also by truncating rigid cusps. (This component $Q'$ is what we refer to in $(5)$ of Theorem \ref{T:ImmersedTurnover}.) In particular, $N$ is hyperbolic. This proves the theorem.  \hfill  \fbox{\ref{T:RegNbhdIsHyp}}\\

\begin{corollary}\label{C:NoNonTOs}
Every cusp neighborhood of $Q$ which is cut out by $f(\cal{T})$ must be rigid. If $\cal{Q}$ is a component of $Q - f(\cal{T})$ with non-elementary fundamental group, then $CC(\cal{Q})$ is either a totally geodesic hyperbolic turnover, a totally geodesic $1$--sided triangle with mirrored sides, or a $3$--dimensional suborbifold with boundary a collection of hyperbolic turnovers which separate $CC(\cal{Q})$ from $f(\cal{T})$.
\end{corollary} 

\section{Bounding the Type of an Elliptic Immersion}\label{S:GehringMartin}

In this section, we will introduce the tools necessary to prove items $(3)$ and $(4)$ of Theorem \ref{T:ImmersedTurnover}, which refer to the orders of the cone points or dihedral angles of the collection of splitting totally geodesic orbifolds which form the boundary of the ``turnover core'' $Q'$. The collection of tools we need comes from the work of Gehring and Martin \cite{GehringMartin94-1}, \cite{Martin96-1}.

If $l_{1}$ and $l_{2}$ are geodesics in $\bb{H}^{3}$, then denote by $\rho(l_{1},l_{2})$ the length of the unique common perpendicular segment between $l_{1}$ and $l_{2}$, or zero if $l_{1}$ and $l_{2}$ intersect in $\bb{H}^{3} \cup S^{2}_{\infty}$. Let $n$ and $m$ be positive integers with $n \geq \max \{ 3,m \}$. Define $c(n,m)$ by
\begin{equation}\label{E:c(n,m)}
	c(n,m)=
	\begin{cases}
		\sqrt{2\cos(2\p/n) - 1} / 2			&\text{if $n \geq 7$,}\\
		\cos(\p/m) / 2					&\text{if $n=6$ and $m \geq 3$,}\\
		1/\sqrt{8}						&\text{if $n=6$ and $m=2$,}\\
		\sqrt{(\sqrt{5}-1)/16}				&\text{if $n=5$,}\\
		\sqrt{(\sqrt{3}-1)/8}				&\text{if $n=4$,}\\
		\sqrt{(\sqrt{5}-2)/8}				&\text{if $n=3$.}\\
	\end{cases}
\end{equation}
Now let $g$ and $h$ be elliptic isometries of order $n$ and $m$, respectively, which generate a Kleinian group. Gehring and Martin show that either $\rho(\text{axis}(g),\text{axis}(h))=0$ or that the following inequality holds \cite[Theorem 6.19]{GehringMartin94-1}:
\smallskip
\begin{equation}\label{E:d(n,m)}
	\rho(\text{axis}(g),\text{axis}(h)) \geq \sinh^{-1} \left( {c(n,m) \over \sin(\p/n) \sin(\p/m)} \right).
\end{equation}
\smallskip
Denote the right-hand side of (\ref{E:d(n,m)}) by $\delta(n,m)$. Observe that it provides a lower bound for the ``axial distance'' between any two elliptic types, i.e., whenever the axes of any two elliptic isometries of orders $n$ and $m$ come closer than $\delta(n,m)$, it is necessary that the group generated by these isometries be elementary. When $n \geq m \geq 7$, any two elliptic elements of orders $n$ and $m$ which generate a discrete group cannot have axes which meet in $\bb{H}^{3} \cup S_{\infty}^{2}$, so $\delta(n,m)>0$ always gives a lower bound for axial distance in this case. We will be most interested in the case $n=m \geq 7$, for which we have the strictly increasing function \cite[Example 8.11]{GehringMartin94-1}, \cite[Theorem 2.2]{Martin96-1}:
\smallskip
\begin{equation}\label{E:d(n)func}
	\delta(n,n) = 2 \cosh^{-1} \left( {1 \over 2 \sin (\p/n)} \right).
\end{equation}
\smallskip
We will show in Section \ref{S:Calcs} that the injectivity radius at any point of a turnover is bounded above by the maximal radius $r_{max}$ of an embedded disk in a hyperbolic thrice-punctured sphere (Proposition \ref{P:InjBound}). We also calculate the value of $r_{max}$:
\smallskip
\begin{equation}
	r_{max} = \ln \left( {2 + \sqrt{7} \over \sqrt{3}} \right) = 0.986647...
\end{equation}
\smallskip

Recall that $\Pi_{T}$ denotes the unique plane stabilized by the turnover subgroup $T = T(p,q,r) \cong \p_{1}(\cal{T}) \leq \C \cong \p_{1}(Q)$, and $\Phi\co \bb{H}^{3} \to Q$ the covering projection. Suppose that $\gamma \in \C$ is an elliptic element of order $n \geq 7$, with axis$(\gamma) \cap \Pi_{T} \neq \emptyset$. Then $\text{axis}(\gamma)$ and $\Pi_{T}$ are either perpendicular or not. Suppose they are not perpendicular. We cannot have $\text{axis}(\gamma) \subset \Pi_{T}$, because then some $T$--translate of $\text{axis}(\gamma)$ would intersect $\text{axis}(\gamma)$, and this would violate discreteness because $\gamma$ has order at least $7$. So the intersection is a point (corresponding to a smooth point of $\cal{T}$, but a singular point of $Q$), and the action of $\gamma$ provides the local model for part of the immersion $f\co \cal{T} \to Q$ (i.e., $f(\cal{T}) \subset Q$ near $\Phi(\text{axis}(\gamma) \cap \Pi_{T})$ looks like a plane ``spun'' around an oblique line). Since $n \geq 7$, we must have that $\rho(\text{axis}(\gamma), \text{axis}(\alpha \gamma)) \geq \delta(n,n)$ for all $\alpha \in \C - \Stab(\text{axis}(\gamma))$, that is, $\text{axis}(\gamma)$ has a tubular neighborhood of radius $\delta(n,n)/2$ which is disjoint from any translate by a deck transformation which does not stabilize it. This is equivalent to $\Phi(\text{axis}(\gamma))$ having an embedded (orbifold) tubular neighborhood of radius $\delta(n,n)/2$ in $Q$. But once we have $\delta(n,n)/2 > r_{max}$, there will be an element $\alpha \in T$ for which the tubular neighborhoods of radius $\delta(n,n)/2$ around $\text{axis}(\gamma)$ and $\text{axis}(\alpha \gamma)$ will not be disjoint. We therefore conclude that $\delta(n,n)/2$ must be less than $r_{max}$, and an easy calculation shows that this occurs only for $n \leq 9$.

Now suppose that $\text{axis}(\gamma)$ meets $\Pi_{T}$ in a right angle. Then the group $\langle \gamma, T \rangle$ stabilizes $\Pi_{T}$, and so is an orientation-preserving Fuchsian group. It is known that turnover subgroups are maximal among orientation-preserving Fuchsian groups, so in fact $\langle\gamma, T\rangle$ is a turnover subgroup of $\pi_{1}(Q)$. There are two ways in which this can occur. One is that $\langle \gamma, T \rangle = T,$ in which case we have $\gamma \in T$, and so $\text{axis}(\gamma)$ is contained in the axis of a (maximal) elliptic element of order $p,q,$ or $r$. The other way that $\langle \gamma, T \rangle$ can be a turnover subgroup is that $T$ is contained in some turnover \emph{super}subgroup $T'(l,m,n) \leq \p_{1}(Q)$, or equivalently, when $f(\cal{T})$ covers a smaller immersed turnover $f'\co \cal{T}' \to Q$ (recall that the hypothesis of the main theorem is that $f(\cal{T})$ does not cover an \emph{embedded} turnover). If this happens, then we must know what types of turnovers can be covered by $\cal{T}$. The subgroups and supergroups of a turnover group can be determined from Table \ref{Ta:TurnoverSupergroups}, whose data is collected from Singerman \cite{Singerman72}. 
\begin{table}
  \begin{center}
    \begin{tabular}{ | c | c | c | c | }
      \hline
      \rule[-8pt]{0pt}{22pt}
      $T(l,m,n)$	&	$\geq T(p,q,r)$	&	Index	&	Normal	\\ \hline
      $(3,3,t)$	&	$(t,t,t)$		&	$3$		&	Yes		\\ \hline
      $(2,3,2t)$	&	$(t,t,t)$		&	$6$		&	Yes		\\ \hline
      $(2,s,2t)$	&	$(s,s,t)$		&	$2$		&	Yes		\\ \hline
      $(2,3,7)$	&	$(7,7,7)$		&	$24$		&	No		\\ \hline
      $(2,3,7)$	&	$(2,7,7)$		&	$9$		&	No		\\ \hline
      $(2,3,7)$	&	$(3,3,7)$		&	$8$		&	No		\\ \hline
      $(2,3,8)$	&	$(4,8,8)$		&	$12$		&	No		\\ \hline
      $(2,3,8)$	&	$(3,8,8)$		&	$10$		&	No		\\ \hline
      $(2,3,9)$	&	$(9,9,9)$		&	$12$		&	No		\\ \hline
      $(2,4,5)$	&	$(4,4,5)$		&	$6$		&	No		\\ \hline
      $(2,3,4t)$	&	$(t,4t,4t)$		&	$6$		&	No		\\ \hline
      $(2,4,2t)$	&	$(t,2t,2t)$		&	$4$		&	No		\\ \hline
      $(2,3,3t)$	&	$(3,t,3t)$		&	$4$		&	No		\\ \hline
      $(2,3,2t)$	&	$(2,t,2t)$		&	$3$		&	No		\\ \hline
    \end{tabular}
    \medskip
    \caption{Turnover Supergroups and Subgroups}\label{Ta:TurnoverSupergroups}
  \end{center}
\end{table}
The data in the first column of the table gives the turnover groups that contain turnover subgroups.  The second column gives the turnover subgroups so contained. Any turnover group not listed in the second column is maximal. By analyzing the table, it is readily seen that the order of $\gamma$ must live in the set 
	\begin{equation}\label{E:AlmostAllConePts}
		\left\{ 2,3,4,5,7,8,9,p,q,r,2p,2q,2r \right\}.
	\end{equation}

Now if $\cal{T}_{i}$ is a turnover or triangle in the splitting collection for $Q'$ as in the main theorem, then the cone points or dihedral angles of $\cal{T}_{i}$ arise from axes of elliptic elements such as $\gamma$ above, that is, elliptic elements whose axes intersect $\Pi_{T}$ in a single point. This is because such an elliptic axis projects directly to a cone point (or dihedral point) in $\partial Q'$. If the cone points of the immersed turnover have orders $p,q,$ and $r$, then by our discussion above, it is possible for the orders of the singular points of any embedded turnovers or triangles in the complement of $f(\cal{T})$ to take values in the set 
	\begin{equation}\label{E:AllConePts}
		\left\{ 2,3,...,9,p,q,r,2p,2q,2r \right\}.
	\end{equation}
Items $(3)$ and $(4)$ of the main theorem follow.

\section{Proof of the Main Theorem}\label{S:MainThmProof}

In this section, we will put everything together and prove the main theorem. Let $f\co \cal{T} \to Q$ be a totally geodesic immersion of a hyperbolic turnover $\cal{T} = \cal{T}(p,q,r)$ in a complete, orientable, hyperbolic $3$--orbifold $Q$. Recall that $p,q,r \in \bb{Z}$ satisfy $2 \leq p \leq q \leq r$ and $\frac{1}{p} + \frac{1}{q} + \frac{1}{r} < 1$. Assume, as in Theorem \ref{T:ImmersedTurnover}, that $f\co \cal{T} \to f(\cal{T})$ is not a covering of an embedded turnover or triangle in $Q$. Let $\{ \cal{T}_{i} \}$ be the collection of hyperbolic turnovers and triangles obtained from the convex cores of components of $Q - f(\cal{T})$ with non-elementary fundamental group. This collection is embedded, pairwise disjoint, and each element of the collection is disjoint from $f(\cal{T})$, which is item $(1)$ of the theorem. Let $S_{j}$ denote a component of  $Q - f(\cal{T})$. Recall that projection from $\partial \overline{S}_{j}$ onto its image in $CC(S_{j})$ strictly decreases area. We have the bound 
	\begin{equation}\label{E:AreaBoundForComplement}
		\bigcup_{j} \Area \left(\partial \overline{S}_{j} \right) \leq 2 \Area \left(\cal{T} \right) 
		= 4\p \left(1-\left(\frac{1}{p} + \frac{1}{q} + \frac{1}{r} \right) \right)
	\end{equation}
on the total area of the boundary of these complementary components.  The reason we do not have equality above is that it may be that $f(\cal{T})$ covers a smaller immersed orbifold in $Q$ (equivalently, $\p_{1}(\cal{T})$ is contained in some larger Fuchsian subgroup).  It follows that the upper bound that we can give on the surface area of the complement of the turnover can be improved the more we know about the Fuchsian groups $G$ such that $\p_{1}(\cal{T}) \leq G \leq \p_{1}(Q)$. These remarks will be useful in Section \ref{S:Apps}.

When $\p_{1}(S_{j})$ is non-elementary, we also have a lower bound of ${\p \over 21}$ for the area of the image of the projection of $\partial \overline{S}_{j}$ onto $\partial CC(S_{j})$, coming from twice the area of the $(2,3,7)$ mirrored triangle, which is the minimal area hyperbolic $2$--orbifold. These two observations imply not only that there is  a \emph{global} upper bound on the number of elements in $\{ \cal{T}_{i} \}$, but also an upper bound (which is at least as strong as the global bound)  on this number  depending only on $p,q$, and $r$. This is item $(2)$ of the theorem.

Now let $Q'$ be the component of $Q - \cup_{i} \cal{T}_{i}$ which contains $f(\cal{T})$. The path metric closure of $Q'$ is homeomorphic to the $3$--orbifold $N$ constructed in Section \ref{S:NbhdHyp}. It follows that the metric closure of $Q'$ is a small hyperbolic orbifold with a totally geodesic boundary component for every element of $\{ \cal{T}_{i} \}$. In Section \ref{S:TurnoverComp}, we showed that we have upper bounds  on the volume of each piece of $Q' - f(\cal{T})$ in terms of the area of $f(\cal{T})$.  We will use the volume estimates (and other claims) provided by Lemmata \ref{L:LimSet0}, \ref{L:LimSet1Rigid}, \ref{L:LimSetInf}, and \ref{L:LimSet2} in the calculation below, depending on the dimension of the convex core of a component $S_{j} \subset Q - f(\cal{T})$, which can be 3, 2, 0 or $-1$ (the latter occurs in the parabolic case; Lemma \ref{L:LimSet2} tells us that $\text{dim } CC(S_{j}) \neq 1$). Recall that $H=1.199678...$ is the positive solution of the equation $x = \coth x.$  We have 
	\begin{align}\label{E:MainTheoremCalc} 
	\Vol(Q') &= \sum_{j} \Vol(S_{j} - CC(S_{j})) \\ \notag \smallskip
	&= \sum_{\text{dim } CC(S_{j}) < 1} \Vol(S_{j}) + 
		\sum_{\text{dim } CC(S_{j}) > 1} \Vol(S_{j} - CC(S_{j})) \\ \notag  \smallskip
	&< \sum_{\text{dim } CC(S_{j}) < 1} {1 \over 2} \cdot \Area(\partial \overline{S}_{j}) + 
		\sum_{\text{dim } CC(S_{j}) > 1} {H \over 2} \cdot \Area(\partial \overline{S}_{j}) \\ \notag \smallskip
	&< {H \over 2} \cdot \sum_{j} \Area(\partial \overline{S}_{j}) \leq {H \over 2} \cdot 
		2 \cdot \Area(\cal{T}) = H \cdot \Area(\cal{T}).
	\end{align}
This is almost all of item $(5)$ of the theorem. The last statement follows from the fact that we have no terms above with $H$ when there are no boundary components for the orbifold $Q'$. This completes the proof of $(5)$ in Theorem \ref{T:ImmersedTurnover}.

We still need to prove the finiteness result, that is, that for a given $p,q,$ and $r$, there are only finitely many possibilities for $Q'$. We use J\o rgensen's Theorem \cite[Theorem 5.12.1]{Thurstonnotes}, \cite[Theorem 5.5]{DunbarMeyerhoff94-1}: Given an upper bound $K$ for volume, there is a finite collection $\cal{O} = \{ O_{1}, O_{2}, ..., O_{k} \}$ of hyperbolic $3$--orbifolds such that any hyperbolic $3$--orbifold with volume $< K$ is obtained from some element of $\cal{O}$ by hyperbolic Dehn surgery. We can apply this theorem to the double $DQ'$ of $Q'$ along its boundary.

Since infinitely many hyperbolic $3$--orbifolds are obtained by hyperbolic Dehn surgery on the elements of $\cal{O}$, we need to prove that only finitely many of these can contain an immersed $(p,q,r)$ turnover. To do this, it is sufficient to show that no surgery with ``large'' coefficients on any element of $\cal{O}$ can yield an orbifold containing an immersed $(p,q,r)$ turnover. 

An embedded neighborhood of a non-rigid cusp in a hyperbolic $3$--orbifold is homeomorphic to either a solid torus minus its core curve or an open solid pillow minus its core  singular curve. Hyperbolic Dehn surgery on a cusped hyperbolic $3$--orbifold is performed by removing an embedded open neighborhood of a cusp, choosing an isotopy class of a closed curve $c$ in the resulting boundary component $V$, and gluing either a solid torus or solid pillow (depending on the cusp) to $V$ so that $c$ is attached to the unique isotopy class of a meridian curve in the solid torus or pillow.  The surgery has ``large'' coefficients if the isotopy class of $c$ has slope ${s \over t}$ (in the universal cover of $V$) and if  $|s| + |t|$ is large.  If hyperbolic Dehn surgery with large coefficients is performed on a cusped hyperbolic $3$--orbifold, then the resulting core curve of the filled solid torus or pillow is either a short orbifold geodesic or an orbifold geodesic with a very small cone angle. Thurston's Hyperbolic Dehn Surgery Theorem (e.g. \cite[Theorem 5.3]{DunbarMeyerhoff94-1}) says that if one excludes a finite number of slopes from each cusp in a cusped hyperbolic $3$--orbifold, then the hyperbolic Dehn surgeries on the remaining slopes all yield hyperbolic $3$--orbifolds. 

Suppose that $DQ'$ is obtained from the collection $\cal{O}$ by hyperbolic Dehn surgery, and that $DQ'$ contains (two copies of) the immersed $(p,q,r)$ turnover $f(\cal{T})$. Let $l \subset DQ'$ be the orbifold geodesic core of a solid torus or pillow coming from a filled cusp. Note that $l$ corresponds to the axis $A$ of a hyperbolic isometry of $\bb{H}^{3}$ which acts on a tubular neighborhood of $A$ by either a translation and rotation or by a dihedral group and rotation.  A result due to Meyerhoff \cite{Meyerhoff87-1} implies that $l$ must have an embedded tubular neighborhood in $DQ'$ and, furthermore, that the radius of such a neighborhood goes to infinity if either the length of $l$ goes to zero or the cone angle on $l$ goes to zero.  Hence, a large hyperbolic Dehn surgery on an element of $\cal{O}$ must yield a $3$--orbifold with an embedded tube of large radius.

We have analyzed the components of $Q' - f(\cal{T})$ in detail, and seen that no such component  supports a simple geodesic. So we must have $l \cap f(\cal{T}) \neq \emptyset$ or $l \cap S \neq \emptyset$, where $S$ is an embedded turnover in $DQ'$ corresponding to a boundary component of $Q'$. Suppose first that $l$ meets such an embedded turnover $S$, and suppose that $l$ is contained in the singular locus of $DQ'$. Then the cone angle at $l$ is $2\pi$ divided by an integer from the set $\{ 2,3,...,9,p,q,r,2p,2q,2r \}$, by item (3) of the main theorem. As a consequence, the cone angle of such a filled cusp cannot be arbitrarily small. The other possibilities are either a non-singular geodesic which meets such an embedded turnover $S$, or else an orbifold geodesic for which $l \cap f(\cal{T}) \neq \emptyset$. But in both of these cases, we may apply the Meyerhoff bound mentioned above to conclude that $l$ can be neither too short nor have a small cone angle, because a hyperbolic turnover has an upper bound on its injectivity radius, and a tube around $l$ with radius larger than this maximal injectivity radius will not be embedded in $DQ'$. So $l$ could not have come from a surgery on any element of $\cal{O}$ with large coefficients. This proves the finiteness claim, and the main theorem.

\section{Fine-Tuning the Search for Immersed Turnovers}\label{S:Miyamoto}

In this section, we will generalize results of Miyamoto \cite{Miyamoto94}, with the goal (in Section \ref{S:Apps}) of limiting the types of embedded turnovers and triangles that can occur in the complement of a few specific immersed turnovers $\cal{T}(p,q,r)$. Let $N$ be the hyperbolic $3$--orbifold obtained by path metric completion of the splitting along the collection $\{ \cal{T}_{i} \}$ of embedded turnovers and triangles from the main theorem. Then $\partial N$ is a collection of hyperbolic turnovers, one for each $\cal{T}_{i}$.

We need some terminology. In the projective model of hyperbolic $n$--space $\bb{H}^{n}$, consider a linearly independent set of $n+1$ points which lie either on the sphere at infinity $S^{n-1}_{\infty}$ or outside of the projective ball. If the line segment between each pair of these points intersects the interior of the projective ball, then the $n+1$ points determine a \emph{truncated $n$--simplex}. This is obtained by taking the infinite volume polyhedron in $\bb{H}^{n}$ spanned by the points, and cutting off the infinite volume ends by the hyperplanes which are dual to the super-ideal vertices. A truncated simplex is \emph{regular} if every edge between two of these truncating planes has the same length. If $r \geq 0$, then define $\rho_{n}(r)$ to be the ratio of the volume of a regular truncated $n$--simplex with edge length $2r$ to the $(n-1)$--volume of its truncated faces. The case of interest to us is dimension three, and we will write $T_{\th}$ for a regular truncated $3$--simplex whose non-truncated faces meet in the angle $\th$, noting that $\th$ and the edge length $2r$ of $T_{\th}$ are related by 
	\begin{equation}\label{E:TruncSideVsAngle}
		\cosh 2r = {\cos \th \over 2 \cos \th - 1}.
	\end{equation}
The truncated faces of regular truncated $3$--simplex of angle $\th$ are equilateral hyperbolic triangles with angle $\th$, and so if $T_{\th}$ has edge length $2r$, then a result due to Miyamoto \cite[Proposition 1.1]{Miyamoto94} implies 
	\begin{align}\label{E:rho3}
		\rho_{3}(r) &= {\Vol(T_{\th}) \over 4(\p-3\th)} \\ \notag
		&= {1 \over 4(\p-3\th)} \left(-8\int_{0}^{\p/4} \ln(2\sin u) \, du - 3\int_{0}^{\th} 
			\cosh^{-1} \left( {\cos t \over 2\cos t - 1} \right) \, dt \right).
	\end{align}
A \emph{return path} in a hyperbolic $n$--orbifold $Q$ with boundary is a geodesic segment in $Q$ which meets $\partial Q$ perpendicularly at both of its end points.  

All of Miyamoto's results are given for hyperbolic manifolds, but some of his arguments do not require any sort of assumptions of no torsion. In particular, the following result is true (cf. \cite[Lemma 4.1]{Miyamoto94}): 

\begin{proposition}\label{P:Miyamoto1}
	If a complete hyperbolic $n$--orbifold $Q$ of finite volume with totally geodesic boundary has a lower bound $l \geq 0$ for the length of its return paths, then $$\Vol(Q) \geq \rho_{n}\left({l \over 2}\right) \Vol(\partial Q).$$
\end{proposition}

Call a return path in $Q$ \emph{closed} if it begins and ends at the same point of $\partial Q$. Such a path must lift to the universal cover of $Q$ as a geodesic segment which meets perpendicularly the axis of an order $2$ elliptic element in $\pi_{1}(Q)$.  With a little bit of careful analysis, we will prove the following result (cf. \cite[Lemma 5.3]{Miyamoto94}): 

\begin{proposition}\label{P:Miyamoto2}
	Let $Q$ be a complete hyperbolic $3$--orbifold with closed totally geodesic boundary.  Then $Q$ has a shortest return path $\gamma$, and there is a positive integer $k$ such that the length of $\gamma$ is at least the edge length of $T_{\th}$, where
	$$\th = {\p \over 3(1-k\chi(\partial Q))}$$
if $\gamma$ is closed and 
	$$\th= {\p \over 3(1-{k \over 2}\chi(\partial Q))}$$
\smallskip
if $\gamma$ is not closed. In both cases, $k>1$ if and only if $\gamma$ is contained in a singular axis of (maximal) order $k$ in $Q$.
\end{proposition}

We can apply these results to our small hyperbolic $3$--orbifold $N$. For a shortest return path $\gamma$ in $N$ and the corresponding truncated regular $3$--simplex $T_{\th}$ of side length $2r$ of Proposition \ref{P:Miyamoto2}, we have, by Theorem \ref{T:ImmersedTurnover} and Proposition \ref{P:Miyamoto1}, the inequalities
	\begin{equation}\label{E:MiyamotoBound}
		{1 \over 4(\p - 3\th)} \Vol(T_{\th}) \Area(\partial N) =
			\rho_{3} \left({2r \over 2} \right) \Area(\partial N) 
	\end{equation}
		$$\leq \Vol(N) < H \cdot \Area(\cal{T}(p,q,r)).$$
(Recall $H=1.199678...$)  It will turn out that we can use these inequalities to limit the number and kinds of turnovers which may appear in $\partial N$, by showing that, for certain types of turnovers in $\partial N$, the lower bound for $\Vol(N)$ is greater than the upper bound.  This will be done in the next section.  To get there, we must first prove the proposition. \medskip

\emph{Proof of \ref{P:Miyamoto2}.} The proof is identical to Miyamoto's proof, provided we take care to consider torsion. However, its elegance is worthy of replication, and we will find use for some of the ideas in Section \ref{S:Apps}.

$Q$ has a shortest return path because its boundary is closed. Denote such a path by $\gamma$, and its length by $l$. Let $\gamma_{1}$ and $\gamma_{2}$ be lifts of $\gamma$ to the universal cover $\tilde{Q} \subset \bb{H}^{3}$ (which is a convex region bounded by geodesic planes), such that $\gamma_{1}$ and $\gamma_{2}$ meet a common component $P$ of $\partial \tilde{Q}$. The argument now diverges according to whether or not $\gamma_{1}$ and $\gamma_{2}$ are equal. There are two ways in which equality can occur. One is that $\gamma$ is contained in the singular locus of $Q$, and the other is that $\gamma$ passes perpendicularly through a singular axis of order $2$ in $Q$. We will first suppose that $\gamma_{1} \neq \gamma_{2}$, for all choices of lifts of $\gamma$. In particular, $\gamma$ is not closed. Let $P_{1}$ and $P_{2}$ denote the other planes of $\partial \tilde{Q}$ which $\gamma_{1}$ and $\gamma_{2}$ meet, respectively, and let $\gamma'$ denote the unique common perpendicular of $P_{1}$ and $P_{2}$. Then the geodesics $\gamma_{1}, \gamma_{2},$ and $\gamma'$ are coplanar and form the alternating sides of an all-right hyperbolic hexagon. Let $d$ be the distance of the side connecting the end points of $\gamma_{1}$ and $\gamma_{2}$ in $P$, and let $l'$ denote the length of $\gamma'$. Observe that $l'$ is no less than $l$. So the all-right hexagon law of hyperbolic cosines  implies
	\begin{equation}\label{E:dInTermsOfl}
		\cosh d = {\cosh^{2} l + \cosh l' \over \sinh^{2} l} \geq {\cosh l \over \cosh l - 1}.
	\end{equation}
Because the above inequality is valid for all choices of lifts of $\gamma$, there are two disjoint disks of radius $r={1\over2}\cosh^{-1}(\cosh l / (\cosh l - 1))$ contained in $\partial Q$ and centered at the end points of $\gamma$.

We now may apply B\"or\"oczky's estimate for circle packings in spaces of constant curvature \cite{Boroczky78}. This says that, given a circle packing of radius $r$ in a space of constant curvature, the density of each disk in its Dirichlet-Voronoi cell (that is, the collection of points lying nearer to the center of the disk than to any other disk in the packing) is at most the density of three mutually tangent disks of radius $r$ in an equilateral triangle spanned by their centers. Let $\th$ denote the angle of a hyperbolic equilateral triangle of side length $2r$. The area of $\partial Q$ is $-2\p \chi(\partial Q)$, and the area of a disk of radius $r$ is $2\p(\cosh r - 1)$, so we have the inequality
	$${4\p(\cosh r - 1) \over -2\p \chi(\partial Q)} \leq {3\th(\cosh r - 1) \over \p - 3\th},$$
which implies that $\th \geq \p/(3(1-\chi(\partial Q)/2))$. By the law of hyperbolic cosines, we obtain
	$$\cosh 2r = {\cos^{2} \th + \cos \th \over \sin^{2} \th} 
	\leq {\cos(\p/(3(1-\chi(\partial Q)/2))) \over 1 - \cos (\p/(3(1-\chi(\partial Q)/2)))}.$$
Now we have 
	$$\cosh l = {\cosh 2r \over \cosh 2r - 1} 
	\geq {\cos (\p/(3(1-\chi(\partial Q)/2))) \over 2 \cos (\p/(3(1-\chi(\partial Q)/2))) - 1}.$$
But equation (\ref{E:TruncSideVsAngle}) implies that $l$ is at least as large as the side length of a regular truncated $3$--simplex with angle $\p/(3(1-\chi(\partial Q)/2))$. This proves the $k=1$ case when $\gamma$ is not closed.

The case when $\gamma$ is closed but not contained in the singular locus follows in precisely the same manner, except that we now take torsion into account. So suppose that $\gamma_{1} = \gamma_{2}$ and that $\gamma$ is not contained in the singular locus of $Q$. Then $\gamma$ must meet an order two axis of $Q$ perpendicularly. If we consider another lift of $\gamma$ which \emph{is} distinct from $\gamma_{1}$, then the same argument applies as above, except that we have only \emph{one} disk packed into $\partial Q$ rather than two. The same calculation yields a bound for the length $l$ of $\gamma$ in terms of $\th$ which is equal to the claim in the closed case with $k=1$.

We are left with the case that $\gamma$ is contained in the singular set of $Q$. Then there is an elliptic isometry in $\p_{1}(Q)$ of some maximal order $k$ whose axis contains $\gamma_{1}$. Suppose first that $\gamma_{1}=\gamma_{2}$ and that $\gamma$ is closed.  As in the last case, we pick a lift of $\gamma$ which is distinct from $\gamma_{1}$, and we conclude that we have one orbifold disk with cone point of order $k$ packed in $\partial Q$. In order to get an estimate on the radius of this orbifold disk, we again appeal to B\"or\"oczky's estimate.  Since our packing in $\partial Q$ is with one orbifold disk, the Dirichlet-Voronoi cell of a disk in the packing obtained by lifting to the universal cover has area equal to $k\Area(\partial Q)$ (because the disk is centered around the lift of an order $k$ cone point).  We therefore obtain, as above, the estimate for $r$ as
	$${2\p(\cosh r - 1) \over k(-2\p \chi(\partial Q))} \leq {3\th(\cosh r - 1) \over \p - 3\th}.$$
We now continue as in the proof of the first case, and the value for the bound on $l$ in terms of $\th$ so obtained is equal to the claim in the closed case with $k \neq 1$.

In the final case, we have two orbifold disks with cone points of order $k$ packed in $\partial Q$ (which is not necessarily connected).  In the universal cover of $\partial Q$, we again have a packing by disks of some radius $r$. Looking at two lifts of the two orbifold disks packed into $\partial Q$, we have two associated Dirichlet-Voronoi cells $V_{1}$ and $V_{2}$, each of whose area is $k$ times the area of  its image in $\partial Q$. We again apply B\"or\"oczky's result to this packing of two disks in two cells and obtain the estimate for $r$ as
	$${4\p(\cosh r - 1) \over k(-2\p \chi(\partial Q))} = {4\p(\cosh r - 1) \over \Area(V_{1} \cup V_{2})} 
	\leq {3\th(\cosh r - 1) \over \p - 3\th},$$
and the same analysis completes the proof.  \hfill \fbox{\ref{P:Miyamoto2}}\\

\section{Applications and Examples}\label{S:Apps}

In this section, we will apply the preceding results to some specific hyperbolic turnovers, as well as provide some examples which show that all of the hypotheses of the main theorem are necessary. In the last subsection, we will use the main theorem to correct an error in a result of Maclachlan \cite{Maclachlan96-1}, and to complete the classification of turnover subgroups in the arithmetic cocompact hyperbolic  tetrahedral reflection groups.

\subsection{Immersed $(2,4,5)$ Turnovers}\label{SS:245TOs}
Let us analyze the case of an immersed $(2,4,5)$ hyperbolic turnover $\cal{T}$ in a hyperbolic orbifold $Q$. Let $N$ be the small orbifold obtained by splitting along embedded turnovers in $Q$, as in the main theorem. Then $\Vol(N) < H \cdot \Area(\cal{T}) = 0.376890...$, where $H=1.199678...$ is given by Theorem \ref{T:ImmersedTurnover}.  Using inequality (\ref{E:d(n,m)}), Table \ref{Ta:TurnoverSupergroups}, and the fact that the hypotenuse of the $(2,4,5)$ hyperbolic triangle represents the greatest distance (approximately $0.842482$, which is \emph{less than} $\delta(n,5)$ for all $n \geq 6$) separating two points on the triangle, it is easily seen that the only possibilities for the orders of cone points in $\partial N$ are $2,3,4,5$.  Because the projection from the immersed turnover onto any boundary pieces strictly decreases area, and because the total two-sided surface area of $\cal{T}$ in $Q$ is $\p/5$, an easy calculation implies that the only possibilities for $\partial N$ are a $(2,4,5)$ turnover or a $(3,3,4)$ turnover.

By considering cases, we can rule out a $(3,3,4)$ turnover boundary by using Proposition \ref{P:Miyamoto2} and the inequalities (\ref{E:MiyamotoBound}). For example, because the $(3,3,4)$ turnover has only one order $4$ cone point, we can conclude that, if the shortest return path in $N$ meets this cone point, then it must be closed. In this case, we would have 
	$$\th =  {\p \over 3(1-4\chi(\partial N))} = {\p \over 4},$$
and we compute the lower bound for $\Vol(N)$ as $0.428850... > H \cdot \Area(\cal{T})$. Therefore, the shortest return path in $N$ could not be contained in this order $4$ axis. All the other cases follow similarly. In fact, we can rule out every possibility for a shortest return path in $N$ \emph{except} in the case that $\partial N$ is a $(2,4,5)$ turnover and the shortest return path is contained in either the order $4$ axis giving the order $4$ cone point in $\partial N$ or the order $5$ axis giving the order $5$ cone point in $\partial N$. We now rule these out by producing better lower bounds for $\Vol(N)$ than those implied by Proposition \ref{P:Miyamoto2}.

\begin{figure}[htbp]
		\begin{center}
		\psfrag{2pi/5}{{\small ${\displaystyle 2\p \over \displaystyle 5}$}}
		\scalebox{1}{\includegraphics{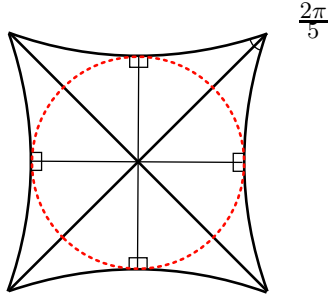}}
		\caption{A lift of the embedded orbifold disk of maximal radius around the order $4$ cone 
			 point in a $(2,4,5)$ hyperbolic turnover}
		\label{F:245TO}
		\end{center}
	\end{figure}

In the order $4$ case, observe that we now know exactly the radius of an embedded orbifold disk around the order $4$ cone point, namely, the length of the segment from the order $4$ cone point to the order $2$ cone point in the $(2,4,5)$ turnover. This distance is given by $\cosh^{-1}(\sqrt{2}\cos(\p/5))$.  See Figure \ref{F:245TO}, which illustrates a lift of this disk to $\bb{H}^{2}$. 
Referring to (\ref{E:dInTermsOfl}) in the proof of Proposition \ref{P:Miyamoto2}, we have that the length $l$ of the shortest return path, if this path is contained in the order $4$ axis, satisfies
	$${\cosh l \over \cosh l - 1} \leq \cosh \left( 2 \cosh^{-1} \left(\sqrt{2}\cos{\p \over 5} \right) \right)
	= 4\cos^{2}{\p \over 5} - 1.$$
This implies
	$$\cosh l \geq {4\cos^{2}\left({\p \over 5} \right) - 1 \over 4\cos^{2}\left({\p \over 5} \right) - 2}.$$
Since we know the relationship of the side length of a regular truncated $3$--simplex to its angle (\ref{E:TruncSideVsAngle}), we can calculate $\th$ so that $l$ is at least the side length of $T_{\th}$. In this case, we have $\th = 0.904556...$, and the corresponding lower bound for $\Vol(N)$ is $0.383986... > H \cdot \Area(\cal{T})$. So the shortest return path can not coincide with an order $4$ cone point in $\partial N$. 

We now consider the order $5$ case. Fix a copy of $\p_{1}(\cal{T})$ and let $\Pi$ be the plane stabilized by it. By assumption, there is an elliptic element of order $5$ which stabilizes a plane $P$ and whose axis meets $\Pi$, and this axis projects to $N$ as a shortest return path.  Applying (\ref{E:d(n,m)}), we have that the axes of two elliptic isometries of order $5$ must be at least $\delta(5,5) \geq 0.736175...$ apart, or else they are equal or they intersect in a point, in order to generate a discrete group. In the last case, the two axes meet in a point which is stabilized by the group of isometries of a regular dodecahedron. Suppose that the axis $A$ of an order $5$ element intersects $\Pi$ (and therefore corresponds to a local part of the immersion of $\cal{T}$), and that the cone angle associated to this element projects to the order $5$ cone point in the $(2,4,5)$ turnover $\partial N$. Then there are several possibilities for the position of $A$. It cannot coincide with the axis of an order $4$ element giving the order $4$ cone point of $\cal{T}$, for then this axis would have to correspond to an element of order at least $20$, and by Table \ref{Ta:TurnoverSupergroups} this cannot occur. So $A$ must be at least $\delta(4,5) \geq 0.626869...$ away from any such order $4$ axis. In this case, we must have that $A$ is close enough to the axis of the order $5$ elliptic giving the order $5$ cone point of $\cal{T}$ so that these axes either intersect or are equal.

\begin{figure}[htbp]
		\begin{center}
		\psfrag{B(N)}{{\small $\partial N$}}
		\psfrag{f(T)}{{\small $f(\cal{T})$}}
		\psfrag{L}{{\small $L$}}
		\psfrag{a}{{\small $a$}}
		\psfrag{b}{{\small $b$}}
		\psfrag{c}{{\small $c$}}
		\scalebox{1}{\includegraphics{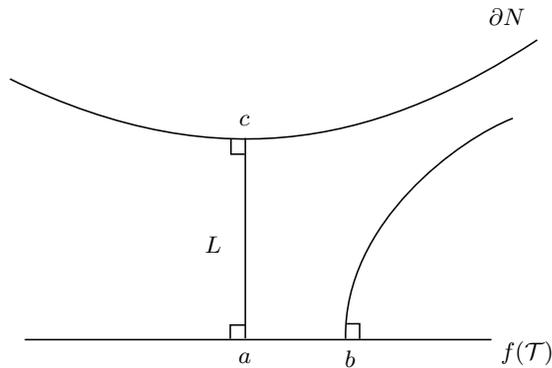}}
		\caption{One way the shortest return path can meet the order $5$ cone point in 
			$\partial N$}
		\label{F:RuleOut5ShortPath1}
		\end{center}
	\end{figure}

Suppose first that they are equal, so that $A$ gives the order $5$ cone points in both $\cal{T}$ and $\partial N$. Then the closed return path of shortest length, if it meets the order $5$ cone point of $\partial N$, must connect these two cone points, one in $\cal{T}$ and one in $\partial N$, and so its length is at least twice the length of this connecting segment (which is the unique common perpendicular for the plane $P$ and the plane $\Pi$ stabilized by our $(2,4,5)$ turnover subgroup), because the return path is closed. Now since the $(2,4,5)$ turnover $\partial N$ is embedded, no other axis giving a cone point of $\cal{T}$ can meet $P$, since two planes have a unique common perpendicular in $\bb{H}^{3}$. This holds in particular for a closest order $4$ axis to $A$. See Figure \ref{F:RuleOut5ShortPath1}, which illustrates the plane containing these two axes. 
The points $a$ and $b$ correspond to order $5$ and order $4$ cone points of $\cal{T}$, respectively,  and $c$ to the order $5$ cone point of $\partial N$. The segment between $a$ and $c$ is contained in the shortest return path, and the ray from $b$ must be disjoint from the geodesic containing $c$. Then by the almost-right quadrilateral law of hyperbolic trigonometry, we have
	$$L = \text{length}(\overline{ac}) > \sinh^{-1} \left( {1 \over \sinh (\text{length}(\overline{ab}))} \right).$$
Since $\overline{ab}$ is the segment between cone points of order $4$ and $5$ in the $(2,4,5)$ turnover, its length is given by $\cosh^{-1}(\cot(\p/5))$, so a calculation gives that $L$ must be bigger than $0.921365...$. Since the length of the shortest return path must be at least twice this, we can calculate the side length for the associated $T_{\th}$, where $\th$ turns out to be $0.938037...$, and the lower bound for $\Vol(N)$ is $0.460222... > H \cdot \Area(\cal{T})$. This proves that the shortest return path cannot connect the order $5$ cone points of $\cal{T}$ and $\partial N$. 

\begin{figure}[htbp]
		\begin{center}
		\psfrag{B(N)}{{\small $\partial N$}}
		\psfrag{A}{{\small $A$}}
		\psfrag{arccos(sqrt(2/(5-sqrt(5))))}{{\small $\cos^{-1} \left( \sqrt{2}/\sqrt{5-\sqrt{5}} \right)$}}
		\scalebox{1}{\includegraphics{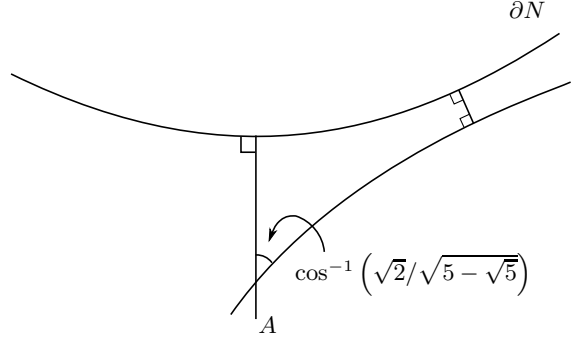}}
		\caption{The other way the shortest return path can meet the order $5$ cone point in 
			$\partial N$}
		\label{F:RuleOut5ShortPath2}
		\end{center}
	\end{figure}

The other possibility is that $A$ intersects in a single point the order $5$ axis giving the cone point of $\cal{T}$.  In this case, this intersection point is stabilized by the orientation-preserving group of symmetries of the dodecahedron, and there is therefore an order $2$ element whose axis meets $A$ in an angle of $\cos^{-1}(\sqrt{2}/\sqrt{5-\sqrt{5}})$.  See Figure \ref{F:RuleOut5ShortPath2}.  
But now we can consider the common perpendicular segment  between $\partial N$ and this order $2$ axis. The path along this common perpendicular which travels from $\partial N$ to this order $2$ and back again will give a shorter return path than the path contained in our order $5$ axis. So the shortest return path could not have been contained in the order $5$ axis.

This analysis implies that there can be no shortest return path in $N$, and so $N$ could have had no boundary. We have therefore shown the following:

\begin{proposition}\label{P:NoEmbTOsIn245Comp}
If a hyperbolic $3$--orbifold $Q$ contains an immersed, non-embedded $(2,4,5)$ turnover, then it can contain no embedded turnovers. In particular, $\Vol(Q) < \Area(\cal{T}(2,4,5)) = \p/10 = 0.314159....$
\end{proposition}

The last claim follows from the last statement in $(5)$ of Theorem \ref{T:ImmersedTurnover}.  The volume bound for such a hyperbolic $3$--orbifold is rather small, and there is such an orbifold \cite{BaskanMacbeath}.  Consider the tetrahedron from the Introduction, shown in Figure \ref{F:BaskMacbIntro}. Recall that $\C_{3}$ denotes the orientation-preserving index two subgroup of the group generated by reflections in the faces of the this tetrahedron, and that the face $ABC$ corresponds to an immersed $(2,4,5)$ turnover in the quotient orbifold $Q_{3} = \bb{H}^{3}/\C_{3}$. The volume of this orbifold is approximately $0.071770$, which is less than a quarter of our volume bound. However, the $(2,4,5)$ turnover subgroup is actually contained in a $\bb{Z}/2\bb{Z}$ extension in $\p_{1}(Q_{3})$ (because it is contained in the group of reflections in the sides of the face $ABC$), and so, by the remarks following (\ref{E:AreaBoundForComplement}), the volume bound of the main theorem can actually be cut in half to $\p/20 = 0.157079...$, which gives a reasonably close approximation for this example. The author conjectures that this is the only example of an immersed, non-embedded $(2,4,5)$ turnover in a hyperbolic $3$--orbifold.

\subsection{Immersed $(2,4,6)$ Turnovers}\label{SS:246TOs}
\begin{figure}[htbp]
		\begin{center}
		\psfrag{A}{{\small $A$}}
		\psfrag{B}{{\small $B$}}
		\psfrag{C}{{\small $C$}}
		\psfrag{D}{{\small $D$}}
		\psfrag{2}{{\small $2$}}
		\psfrag{3}{{\small $3$}}
		\psfrag{4}{{\small $4$}}
		\psfrag{6}{{\small $6$}}
		\psfrag{pi/4}{{\small ${\displaystyle \p \over \displaystyle 4}$}}
		\psfrag{pi/6}{{\small ${\displaystyle \p \over \displaystyle 6}$}}
		\scalebox{1}{\includegraphics{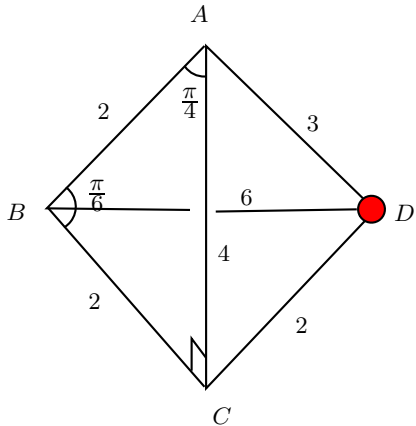}}
		\caption{A noncompact hyperbolic tetrahedron with $(2,4,6)$ triangular face}
		\label{F:246NoncptTet}
		\end{center}
	\end{figure}

Consider the tetrahedron $T_{10}$ in Figure \ref{F:246NoncptTet}.  The notation for the dihedral angles is the same as in the previous example, and it can be similarly realized in $\bb{H}^{3}$ as the noncompact fundamental domain for the discrete group generated by reflections in its sides.  It has one ideal vertex indicated at $D$.
 Let $\C_{10}$ be the orientation-preserving index two subgroup of the reflection group. This is also a Kleinian group, and we have $\Vol(Q_{10}=\bb{H}^{3}/\C_{10}) \approx 0.211446$. Now the face $ABC$ of this tetrahedron is a $(2,4,6)$ triangle, and so $\Gamma_{10}$ contains a $\bb{Z}/2\bb{Z}$ extension of a $(2,4,6)$ turnover group. Clearly $Q_{10}$ can contain no embedded turnovers, and so our volume  bound is $\Area(\cal{T}(2,4,6))/2 = 0.261799...$, which is very close to the volume of $Q_{10}$. We make the following conjecture:

\begin{conjecture}\label{J:NoEmbTOsIn246Comp}
If a hyperbolic $3$--orbifold $Q$ contains an immersed, non-embedded $(2,4,6)$ turnover, then it can contain no embedded turnovers. In particular, $\Vol(Q) < \Area(\cal{T}(2,4,6)) = \p/6 = 0.523598....$
\end{conjecture}

This conjecture is made in the place of calculating, as in the last example, that an immersed $(2,4,6)$ turnover can contain no embedded turnovers in its complement. The author believes that this can be shown, but the case by case analysis necessary is rather complicated. 

\subsection{Immersed $(2,4,p)$ Turnovers for $p \geq 7$}\label{SS:24p>6TOs}
\begin{figure}[htbp]
		\begin{center}
		\psfrag{2}{{\small $2$}}
		\psfrag{T(2,3,764)}{{\small  $\cal{T}(2,3,764)$}}
		\psfrag{embedded}{{\small  embedded}}
		\psfrag{3}{{\small $3$}}
		\psfrag{4}{{\small $4$}}
		\psfrag{764}{{\small $764$}}
		\psfrag{T(2,4,764)}{{\small $\cal{T}(2,4,764)$}}
		\psfrag{immersed}{{\small  immersed}}
		\scalebox{1}{\includegraphics{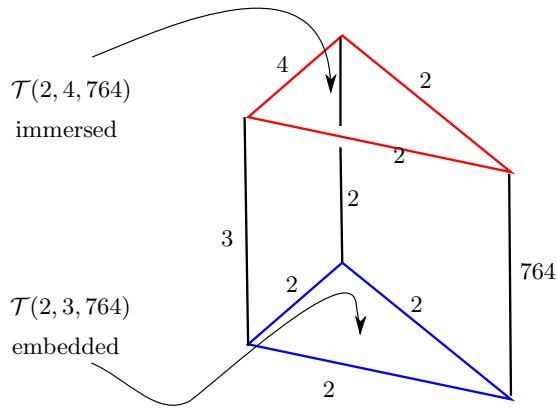}}
		\caption{A hyperbolic prism with immersed and embedded turnovers}
		\label{F:hyperbolicprism}
		\end{center}
	\end{figure}

The conjectural lack of embedded turnovers in the complement of an immersed $(2,4,p)$ turnover does not extend to the case when $p \geq 7$. See Figure \ref{F:hyperbolicprism} (where $p=764 \geq 7$ is arbitrary). 
The index two orientation-preserving subgroup of the group generated by reflections in the sides of this prism is a Kleinian group by Andre'ev's Theorem (e.g. \cite{DunbarHubbardRoeder06}). The corresponding orbifold contains an immersed $(2,4,p)$ turnover (coming from the roof of the prism) and an embedded $(2,3,p)$ turnover in the complement of this immersed turnover (coming from the base). (Recall from the introduction that Martin \cite{Martin96-1} showed that $(2,3,p)$ turnovers are always embedded). Again, the fundamental group of the $3$--orbifold $Q_{2,4,p}$ obtained from this prism contains a $\bb{Z}/2\bb{Z}$ extension of the $(2,4,p)$ turnover group, and so the upper volume bound is ${H \over 2} \Area(\cal{T}(2,4,p))$ ($H=1.199678...$). This bound is best when $p=7$, where we have 
	$$0.325947 \approx \Vol(Q_{2,4,7}) <  {H \over 2} \Area(\cal{T}(2,4,7)) = 0.403810...,$$
and the bound worsens as $p$ tends to infinity, where we have
	$$0.501921 \approx \Vol(Q_{2,4,\infty}) < {H \over 2} \Area(\cal{T}(2,4,\infty)) = 0.942225....$$

We observe that this family of examples demonstrates that the splitting collection of turnovers $\{\cal{T}_{i}\} \subset Q$ from Theorem \ref{T:ImmersedTurnover} is not always empty. The fundamental groups of the orbifolds $Q_{2,4,p}$ in this section are examples of so-called \emph{web} groups.  Their limit sets form web-like Sierpinski gaskets on $S_{\infty}^{2}$ with turnover groups for the component stabilizers.  

\subsection{Turnover Subgroups of Tetrahedral Groups}\label{SS:Maclachlan}
\begin{figure}[htbp]
		\begin{center}
		\psfrag{A}{{\small $A$}}
		\psfrag{B}{{\small $B$}}
		\psfrag{C}{{\small $C$}}
		\psfrag{D}{{\small $D$}}
		\scalebox{1}{\includegraphics{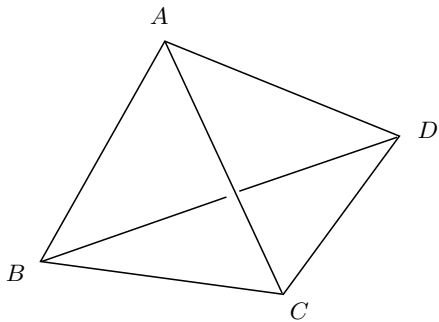}}
		\caption{The tetrahedron $ABCD$}
		\label{F:MaclachlanTemplate}
		\end{center}
	\end{figure}

The principal reference for this subsection is Maclachlan \cite{Maclachlan96-1}, which classifies almost all the turnover subgroups contained in the cocompact hyperbolic tetrahedral groups. There are nine compact hyperbolic tetrahedra such that reflections in the sides of the tetrahedra tile hyperbolic $3$--space \cite[Chapter 7]{Ratcliffe94-1}. We denote them by $T_{i}[l_{1},l_{2},l_{3};m_{1},m_{2},m_{3}]$, where $i$ ranges from $1$ to $9$, and where the $l_{j}$ and $m_{k}$ give the integer submultiples of $\p$ at each edge $AB, BC, AC, CD, AD, BD$ (in this order) for the tetrahedron $ABCD$ given in Figure \ref{F:MaclachlanTemplate}.  
Let $\C_{i}$ denote the orientation-preserving index two subgroup of the group generated by the reflections in the faces of $T_{i}$.  Then the $\C_{i}$ are Kleinian groups which are all  arithmetic except for $\C_{8}$, and Maclachlan classifies most of the turnover subgroups for the eight arithmetic examples by arithmetic means. We will see that it is possible, in principle, to locate these turnover subgroups by purely geometric means. We will focus on $T_{8}[2,3,4;2,3,5]$ and $T_{9}[2,3,5;2,3,5]$.  If we have a turnover subgroup $T(p,q,r) < \C_{i}$, then the plane $\Pi_{T}$ invariant under $T$ will always intersect at least three intersecting elliptic axes in $\bb{H}^{3}$ (corresponding to rotations in three edges of the tetrahedron $T_{i}$). Therefore, at least one of these axes must meet $\Pi_{T}$ non-perpendicularly. It follows that any turnover subgroup $T(p,q,r) < \C_{i}$ will yield an immersed turnover in the quotient $3$--orbifold. 

We begin by correcting an error in Maclachlan's work. He claims to prove \cite[Theorem 5.1]{Maclachlan96-1} that $\cal{O}_{9} = \bb{H}^{3}/\C_{9}$ contains the turnovers $(3,3,5)$ and $(5,5,5)$. He also states that $\cal{O}_{9}$ \emph{may} contain a $(3,5,5)$ turnover.  The volume of $\cal{O}_{9}$ can be computed as approximately $1.004261$. Since $\cal{O}_{9}$ contains no embedded turnovers, the upper bound for its volume given by Theorem \ref{T:ImmersedTurnover} is just the area of any immersed turnover that it contains. But the volume bound in the case of the $(3,3,5)$ turnover is $4\p/15 = 0.837758... < 1.004261$. We conclude that $\cal{O}_{9}$ can contain no $(3,3,5)$ turnover. Maclachlan's claim that $\cal{O}_{9}$ contains a $(5,5,5)$ turnover is based on the fact that $T(5,5,5) < T(3,3,5)$, and so this claim becomes unjustified.

\begin{figure}[htbp]
		\begin{center}
		\psfrag{A}{{\small $A$}}
		\psfrag{B}{{\small $B$}}
		\psfrag{C}{{\small $C$}}
		\psfrag{C'}{{\small $C'$}}
		\psfrag{D}{{\small $D$}}
		\psfrag{2}{{\small $2$}}
		\psfrag{3}{{\small $3$}}
		\psfrag{5}{{\small $5$}}
		\psfrag{2/2}{{\small $2/2$}}
		\psfrag{3/2}{{\small $3/2$}}
		\psfrag{5/2}{{\small $5/2$}}
		\scalebox{1}{\includegraphics{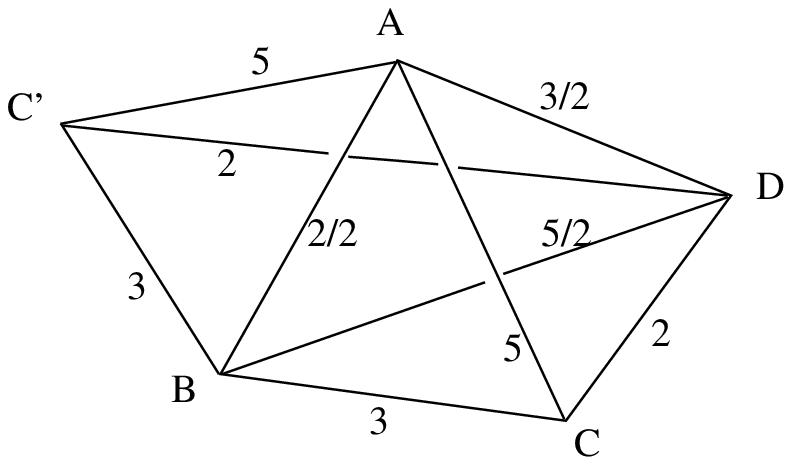}}
		\caption{A fundamental domain for $\C_{9}$ in $\bb{H}^{3}$}
		\label{F:355InT9}
		\end{center}
	\end{figure}

As it turns out, we can prove directly that $\cal{O}_{9}$ contains both a $(3,5,5)$ turnover and a $(5,5,5)$ turnover. We begin with the $(3,5,5)$ turnover. See Figure \ref{F:355InT9}, which illustrates one fundamental domain for $\C_{9}$ in $\bb{H}^{3}$.
The numbers at the edges represent the submultiples of $\p$ of the dihedral angle at that edge  (so, for instance, the edge $BD$ is labeled $5/2$ because the dihedral angle at that edge is $2\p/5$). We remark that this illustration is, obviously, not in true perspective. Note that each edge also represents the axis of an elliptic isometry in $\C_{9}$, whose order is the numerator of the edge label. We will prove the claim by finding three planes whose pairwise lines of intersection are composed of axes of elliptic elements that generate a $(3,5,5)$ turnover subgroup.

\begin{figure}[htbp]
		\begin{center}
		\psfrag{A}{{\small  $A$}}
		\psfrag{B}{{\small  $B$}}
		\psfrag{C}{{\small $C$}}
		\psfrag{C'}{{\small $C'$}}
		\psfrag{D}{{\small $D$}}
		\psfrag{2}{{\small $2$}}
		\psfrag{3}{{\small  $3$}}
		\psfrag{5}{{\small $5$}}
		\psfrag{2/2}{{\small  $2/2$}}
		\psfrag{3/2}{{\small  $3/2$}}
		\psfrag{5/2}{{\small  $5/2$}}
		\psfrag{pi/5}{{\small  $\displaystyle {\p \over 5}$}}
		\psfrag{L1}{{\small  $L_{1}$}}
		\psfrag{L2}{{\small  $L_{2}$}}
		\psfrag{L3}{{\small $L_{3}$}}
		\psfrag{L4}{{\small $L_{4}$}}
		\scalebox{1}{\includegraphics{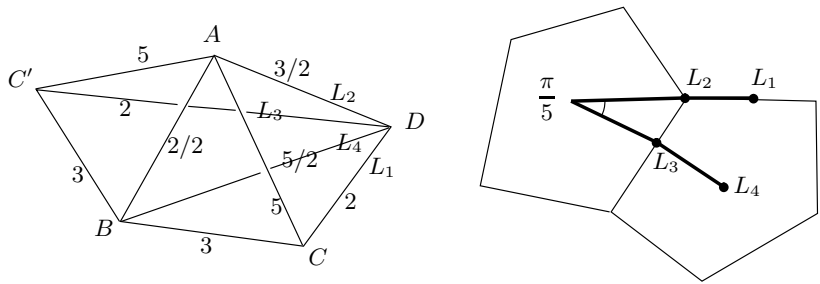}}
		\caption{The intersection of planes in the link of the vertex $D$}
		\label{F:235LinkFor355TO}
		\end{center}
	\end{figure}

Let $\measuredangle(X,Y)$ denote the dihedral angle of the intersecting planes $X$ and $Y$. It is easily seen from the figure that $\measuredangle(ACD, ACC')= \p/5$ and $\measuredangle(ACC',BC'D) = \p/3$, and that these pairs of planes intersect in order $5$ and order $3$ elliptic axes, respectively. The planes $BC'D$ and $ACD$ contain $D$ in their intersection, and we can calculate their line of intersection and the dihedral angle as follows. The point $D$ is stabilized by a $(2,3,5)$ spherical turnover subgroup $G \leq \C_{9}$. The plane $ACD$ is spanned by lines $L_{1} \supset CD$ and $L_{2} \supset AD$, which are the axes  of elements in $G$ of order $2$ and $3$, respectively, and which are as close as possible in terms of their angle of intersection. Similarly,  from the plane $BC'D$ we obtain closest possible axes $L_{3} \supset C'D$ and $L_{4} \supset BD$ of elements of order $2$ and $5$, respectively. To calculate the angle $\measuredangle(ACD,BC'D)$, then, we simply calculate the angle between the planes spanned by these axes in the spherical link of $D$. See Figure \ref{F:235LinkFor355TO}. 
The image at left is just Figure \ref{F:355InT9} with the axes $L_{i}$ labeled, and the image at right shows two faces of the spherical dodecahedron around $D$, corresponding to the $(2,3,5)$ spherical turnover subgroup $G$. The labels $L_{i}$ are meant to indicate the intersection of $L_{i}$ with the link. Now it is easily seen that $\measuredangle(ACD,BC'D) = \p/5$, and that this intersection occurs along the axis of an order $5$ elliptic element. So we have three planes which pairwise intersect in angles giving a $(3,5,5)$ triangle, with elliptic isometries of the appropriate orders at the vertices. These isometries then generate a $(3,5,5)$ turnover subgroup of $\C_{9}$.

\begin{figure}[htbp]
		\begin{center}
		\psfrag{A}{{\small $A$}}
		\psfrag{B}{{\small $B$}}
		\psfrag{C}{{\small $C$}}
		\psfrag{C'}{{\small $C'$}}
		\psfrag{D}{{\small $D$}}
		\psfrag{D'}{{\small $D'$}}
		\psfrag{A'}{{\small $A'$}}
		\psfrag{2}{{\small $2$}}
		\psfrag{3}{{\small $3$}}
		\psfrag{3/3}{{\small $3/3$}}
		\psfrag{5}{{\small $5$}}
		\psfrag{2/2}{{\small $2/2$}}
		\psfrag{3/2}{{\small $3/2$}}
		\psfrag{5/2}{{\small $5/2$}}
		\psfrag{5/3}{{\small $5/3$}}
		\scalebox{1}{\includegraphics{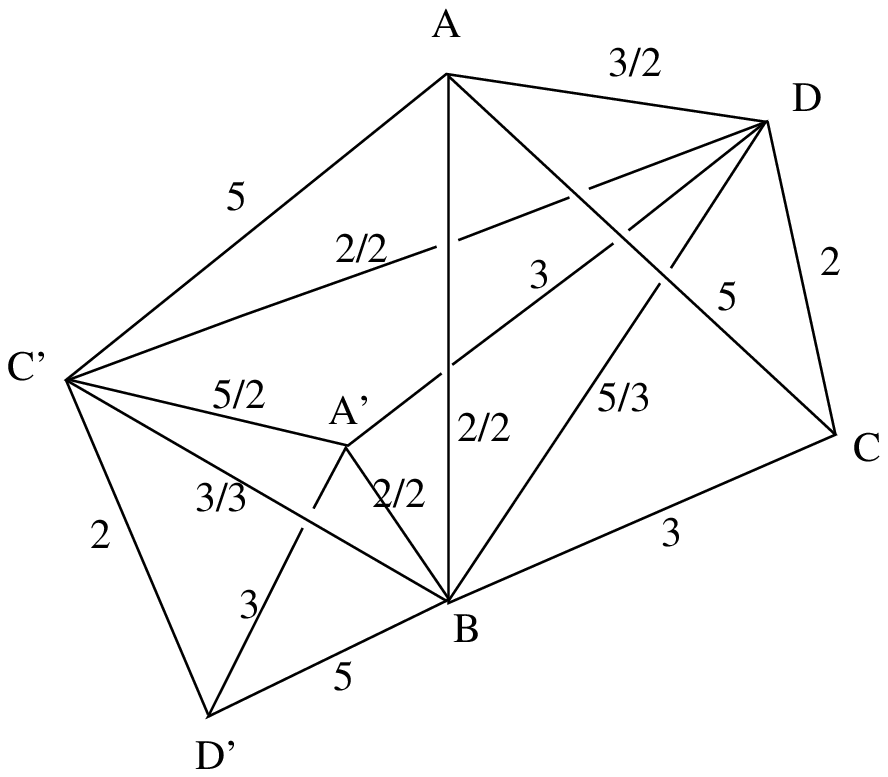}}
		\caption{Two fundamental domains for $\C_{9}$ in $\bb{H}^{3}$}
		\label{F:555InT9}
		\end{center}
	\end{figure}

The same method can also be used to locate a $(5,5,5)$ turnover in $\cal{O}_{9}$. See Figure \ref{F:555InT9}, which illustrates two copies (that is, four tetrahedra) of a fundamental domain for the action of $\C_{9}$.
Since the edge labels indicate the submultiples of $\pi$ for the polyhedron pictured, we observe from the figure that the points $A,B,C,C',$ and $D'$ are coplanar, and similarly for $A',B,D,$ and $D'$. It is then easily seen from the figure that we have $\measuredangle(ACD',ACD) = \measuredangle(ACD',A'DD') = \pi/5$. Additionally, by the method illustrated in Figure \ref{F:235LinkFor355TO}, it is readily seen that $\measuredangle(ACD,A'DD') = \pi/5$. All of these intersections occur along order $5$ axes of elements in $\C_{9}$, and so as in the last case we conclude that $\cal{O}_{9}$ contains a $(5,5,5)$ turnover.

The method just described can be used to verify many of the examples of turnover subgroups in tetrahedral groups given by Maclachlan. It is also possible, in principle, to describe a geometric algorithm to locate all the examples. Because the lengths of the sides of the possible turnovers in these tetrahedral groups are bounded, and because the dihedral edge at each elliptic axis in these tetrahedra is fixed by a translation in the universal cover $\bb{H}^{3}$, the search for axes which are mutually at the appropriate distance apart from one another will be a search over a compact subset of hyperbolic $3$--space. We will not construct an explicit algorithm.  However, we note that the only possibilities for turnover subgroups of $\C_{9}$ are $(2,5,5), (3,3,5), (3,5,5),$ and $(5,5,5)$ (because of the orders of the edges of $T_{9}$), and since we can rule out the first two of these (the (2,5,5) case follows by our volume bound, exactly as in the (3,3,5) case) and prove the existence of the other two, we have the following:

\begin{proposition}\label{P:TOsInT9}
The turnover subgroups of the group $\C_{9}$ are $T(3,5,5)$ and $T(5,5,5)$.
\end{proposition}   

Turning to $T_{8}$, we have $\Vol(\cal{O}_{8}= \bb{H}^{3}/\C_{8}) \approx 0.717306$. The orders of the cone points for any turnover in $\cal{O}_{8}$ must come from the dihedral angles at the edges of $T_{8}$, and so the list of possible turnovers is
    $$(2,4,5), (2,5,5), (3,3,4), (3,3,5), (3,4,4),$$ 
    $$(3,4,5), (3,5,5), (4,4,4), (4,4,5), (4,5,5), (5,5,5).$$ 
The volume bound of the main theorem rules out the left-most three turnovers from the top row. Of those that remain, only $(3,4,5)$ is non-arithmetic. It can be shown that $\cal{O}_{8}$ contains this turnover, using the exact same analysis (and the same figures with slightly different edge labels) as in the $T(3,5,5) < \C_{9}$ case above.

As in the $T(5,5,5) < \C_{9}$ case, one can prove the existence of a $(4,5,5)$ turnover in $\cal{O}_{8}$ as well. Once again, the proof is identical, and the figure required to carry out the proof differs from Figure \ref{F:555InT9} only in some of the edge labels. 

The group $\C_{8}$ may contain other arithmetic turnover subgroups from the list above. However, in the absence of an explicit algorithm (whose existence was described above), we note that the presence of the arithmetic turnover subgroup $T(4,5,5)$ provides interesting information about the orbifold $\cal{O}_{8}$, due to a recent result of Long-Lubotzky-Reid \cite[Proposition 4.1]{LongLubReid07-1}:

\begin{proposition}\label{P:TOsInT8}
The orbifold $\cal{O}_{8}$ cannot contain immersed turnovers of type $(2,4,5)$, $(2,5,5),$ or $(3,3,4)$. It contains immersed $(3,4,5)$ and $(4,5,5)$ turnovers, and the existence of the latter implies that $\cal{O}_{8}$ has a tower of principle congruence covers with Property $\tau$.
\end{proposition}

\section{Calculations}\label{S:Calcs}

In this section, we will prove some results about hyperbolic turnovers and thrice-punctured spheres, as well as determine an interesting new isoperimetric constant. Our first result is a bound on the diameter of a maximally embedded disk in a turnover.

\begin{proposition}\label{P:InjBound}
Let $\cal{T}(p,q,r)$ be a hyperbolic turnover. Then the radius of an embedded disk in $\cal{T}$ is less than  
	\begin{equation}\label{E:3PuncInj}
		\ln \left( {2 + \sqrt{7} \over \sqrt{3}} \right).
	\end{equation} 
\end{proposition}

\emph{Proof of \ref{P:InjBound}.} We will prove the result by proving that the maximal radius $r_{max}$ of an embedded disk in a hyperbolic thrice-punctured sphere is given by (\ref{E:3PuncInj}), and that the maximal radius of an embedded disk in a hyperbolic turnover must be less than $r_{max}$.

\begin{figure}[htbp]
		\begin{center}
		\psfrag{-2}{{\small $-2$}}
		\psfrag{-1}{{\small $-1$}}
		\psfrag{0}{{\small $0$}}
		\psfrag{2}{{\small $2$}}
		\psfrag{g}{{\small $g$}}
		\psfrag{y = sqrt(3)}{{\small $\sqrt{3}$}}
		\scalebox{1}{\includegraphics{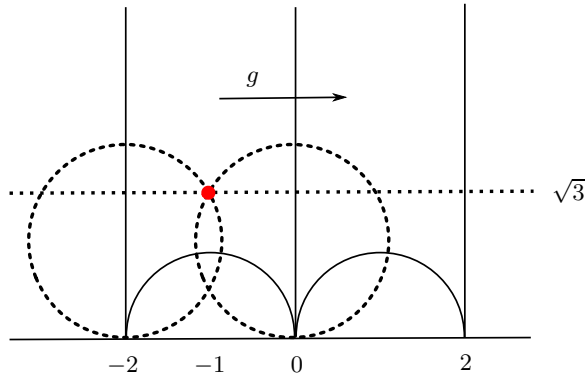}}
		\caption{A fundamental domain for a thrice-punctured sphere, with three 
			isometric horospheres}
		\label{F:3PuncOrtho}
		\end{center}
	\end{figure}

We begin with the thrice-punctured sphere $S$. Consider a fundamental domain in $\bb{H}^{2}$ for $S$, normalized as in Figure \ref{F:3PuncOrtho}. 
We have parabolic isometries which fix $0, -2,$ and $\infty$, the latter acting on the upper half-space as translation by $4$. Denote this translation by $g\co z \mapsto z + 4$. It is an easy computation to show that a rotation $\sigma$ thru $2 \p / 3$ radians centered at the  ``orthocenter'' $(-1,\sqrt{3})$ of the ideal triangle on the left permutes the points $0, -2, \infty$ cyclically. We consider the horoballs based at $0, -2,$ and $\infty$, each passing through $(-1,\sqrt{3})$, as in the figure. These horoball regions are permuted by $\sigma$ with their centers, and any point in $S$ lifts to a point contained in at least one of these horoballs. It is therefore the case that moving in any direction away from $(-1,\sqrt{3})$ is equivalent, up to the action of $\sigma$, to increasing the Euclidean height from $\sqrt{3}$. However, if a point in the upper half-space has coordinates $(a,b)$, then a computation shows that the distance to its translate by $g$ is given by 
	$$2 \ln \left( {2 + \sqrt{b^{2} + 4} \over b} \right),$$
which is less than twice the quantity in (\ref{E:3PuncInj}) when $b> \sqrt{3}$. Any point in the fundamental domain for $S$ is taken to its closest translate by the action of $g$, up to the isometry $\sigma$. So the largest embedded disk in $S$ lifts to a disk based at $(-1,\sqrt{3})$. But, as observed above, a disk based at $(-1,\sqrt{3})$ will meet its closest translates by $\p_{1}(S)$ \emph{exactly} when it has radius given by (\ref{E:3PuncInj}).  

To complete the proof, we construct a contraction from the thrice-punctured sphere to the smooth part of a hyperbolic turnover. Consider a triangle $\triangle$ in the projective model of $\mathbb{H}^{2}$, positioned so that its circumcenter lies at the center of the disk. Then the vertices of $\triangle$ lie on a circle of that is concentric with the unit disk, and $\triangle$ is the (Euclidean) homothetic image of an ideal triangle. It is easily seen that the Euclidean homothety taking this ideal triangle to $\triangle$ is a contraction for the hyperbolic metric. If we double along the boundary of $\triangle$, then this gives a contraction from the thrice-punctured sphere to the double of $\triangle$ that restricts to a homeomorphism on the smooth part of the double of $\triangle$.

Now if the maximal radius of an embedded ball in a hyperbolic turnover is greater than or equal to $r_{max}$, then we have an immediate contradiction. For consider the center $p$ of such a ball, and let $\ell$ be the shortest loop through $p$. Then the length of $\ell$ is at least $2r_{max}$. But then pulling $\ell$ back to the thrice-punctured sphere, decreasing the length of the resulting loop (which we can do because the thrice-punctured sphere has maximal injectivity radius $r_{max}$), and contracting back to the turnover will produce a shorter path through $p$. This proves the proposition. \hfill \fbox{\ref{P:InjBound}}\\

\begin{figure}[htbp]
		\begin{center}
		\psfrag{H}{{\small $H$}}
		\psfrag{2}{{\small $2$}}
		\psfrag{1.66711...}{{\small $1.667113...$}}
		\psfrag{AoverV}{{\small $\displaystyle {\Area(S(H)) \over \Vol(B(H))}$}}
		\scalebox{1}{\includegraphics{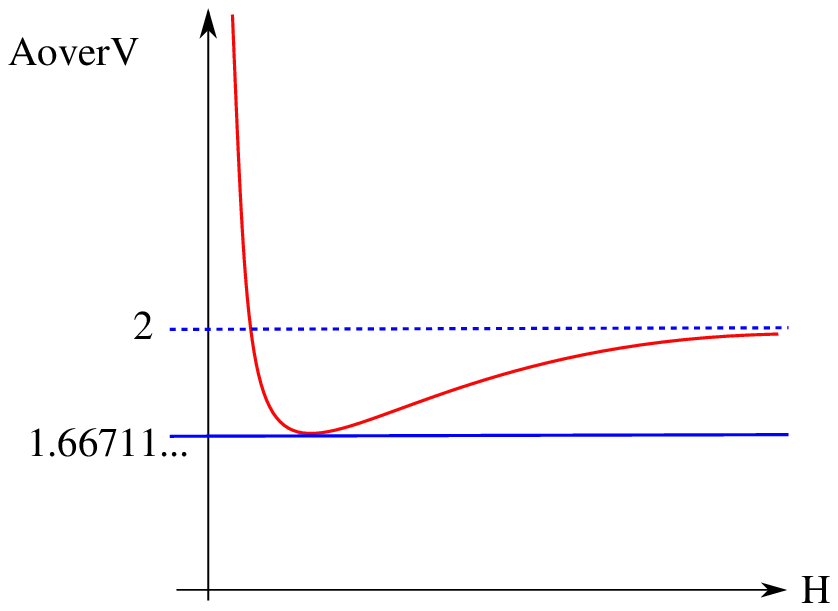}}
		\caption{The function $H \mapsto \displaystyle {4 \cosh^{2} H \over \sinh 2H + 2H }$}
		\label{F:GraphOfArea/Vol}
		\end{center}
	\end{figure}

Our final calculation provides the volume bound of a nice room in terms of the area of its ceiling, which is used in Corollary \ref{C:RoomVolBound}. In the terminology of Section \ref{S:Isop}, let $F$ be a floor in a hyperbolic plane $\Pi \subset \bb{H}^{3}$, and let $B$ and $S$ be a nice room over $F$ and its ceiling, respectively. Let $H$ be the height of $B$.  Recall that the metric for $\bb{H}^{3}$ using Fermi coordinates based on $\Pi$ is given by $dh^{2} + \cosh^{2} h (dr^{2} + \sinh^{2} r d\th^{2})$, where $r$ and $\th$ describe polar coordinates on $\Pi$ and $h$ is perpendicular distance to $\Pi$. Write $dA = \sinh^{2}r \, drd\th$ for the area form on $\Pi$. We have
	\begin{align}\label{E:Area(S)/Vol(B)}
		{\Area(S) \over \Vol(B)} &= { \iint\limits_{F} \cosh^{2}H \, dA  \over 
			\iint\limits_{F} \int_{0}^{H} \cosh^{2} h \, dh \, dA } \\ \notag
		&= {4 \cosh^{2} H \over \sinh 2H + 2H }.
	\end{align}
The function of $H$ given above has interesting properties. It limits to infinity as $H$ approaches zero, and to $2$ as $H$ approaches infinity. It is \emph{not} monotonic, however. It has a minimum value of  $2/H = 1.667113...$ when $H$ is the positive solution of $\coth x = x$, and it decreases from infinity to this value and then increases to $2$ from this value. As a side remark, the author would be very interested in understanding the significance of the critical point of this function, as it seems to defy intuition for this ratio to be non-monotone. Figure \ref{F:GraphOfArea/Vol} shows a graph of this function. In particular, we have
	$${\Area(S) \over \Vol(B)} \geq 1.667113...,$$
which is the result necessary for Corollary \ref{C:RoomVolBound}.

\bibliographystyle{hyperamsplain}
\bibliography{refs}

\end{document}